\documentclass[reqno,10pt]{article}
\usepackage{amsmath, amsthm, graphicx,amsfonts, amssymb,color}
\usepackage{bm}
\usepackage{hyperref}
\RequirePackage{amsthm,amsmath,amsfonts,amssymb,mathtools}
\RequirePackage[numbers]{natbib}

\usepackage{bm}
\usepackage{hyperref}
\usepackage{mathrsfs}

\newtheorem{thm}{Theorem}[section]
\newtheorem{lemma}[thm]{Lemma}

\newtheorem{corol}[thm]{Corollary}
\newtheorem{propos}[thm]{Proposition}
\newtheorem{rema}[thm]{Remark}
\newtheorem{defini}[thm]{Definition}
\newtheorem{example}[thm]{Example}

\def\bex{\begin{example}}
\def\eex{\end{example}}
\def\bd{\begin{defini}}
\def\ed{\end{defini}}
\def\bp{\begin{propos}}
\def\ep{\end{propos}}
\def\bt{\begin{thm}}
\def\et{\end{thm}}
\def\bco{\begin{corol}}
\def\eco{\end{corol}}
\def\bl{\begin{lemma}}
\def\el{\end{lemma}}
\def\br{\begin{rema}}
\def\er{\end{rema}}
\def\be{\begin{equation}}
\def\ee{\end{equation}}
\def\ba{\begin{array}}
\def\ea{\end{array}}
\def\bena{\begin{eqnarray}}
\def\eena{\end{eqnarray}}

\def\P{{\mathbb P}}

\def\R{{\mathbb R}}
\def\Z{{\mathbb Z}}

\def\ln{\log}
\def\1{I}
\def\imath{\textbf{i}}
\def\jmath{\textbf{j}}
\def\chi{\zeta}
\def\fA{{\cal A}}

\def\hS{{\mathscr S}}

\def\hC{{\mathscr C}}

\def\a{{\alpha}}

\def\QED{\hfill$\square$\vskip 3mm}

\def\Dp{\displaystyle}
\def\Df{\Dp\frac}

\def\({\left(}
\def\){\right)}

\begin{document}


\title{\huge An Explicit Description of Extreme Points of the Set of Couplings with Given Marginals: with Application to Minimum-Entropy Coupling Problems
\\[5mm]
\footnotetext{*Correspondence author}
\footnotetext{Research supported in part by the Natural Science Foundation of China (under
grants 11471222, 61973015)}}
\author{Ya-Jing Ma$^1$,\  Feng Wang$^1$,\  Xian-Yuan Wu$^{1*}$,\  Kai-Yuan Cai$^2$}\vskip 10mm

\date{}
\maketitle

\begin{center}
\begin{minipage}{13cm}
{\footnotesize
$^1$School of Mathematical Sciences, Capital
Normal University, Beijing, 100048, China. Emails:
\texttt{mayajing121@126.com},\ \texttt{wangf@cnu.edu.cn},\ \texttt{wuxy@cnu.edu.cn}

$^2$Department of Automatic Control, Beijing University of Aeronautics and Astronautics, Beijing, 100191, China. Email:\texttt{kycai@buaa.edu.cn}
}
\end{minipage}
\end{center}
\vskip 5mm
\begin{center} \begin{minipage}{14.6cm}
{\small {\bf Abstract}: \it Given probability distributions ${\bf p}=(p_1,p_2,\ldots,p_m)$ and ${\bf q}=(q_1,q_2,\ldots, q_n)$ with $m,n\geq 2$, denote by ${\cal C}(\bf p,q)$ the set of all couplings of $\bf p,q$, a convex subset of $\R^{mn}$. Denote by ${\cal C}_e({\bf p},{\bf q})$ the finite set of all extreme points of ${\cal C}(\bf p,q)$. It is well known that, as a strictly concave function, the Shannan entropy $H$ on ${\cal C}(\bf p,q)$ takes its minimal value in ${\cal C}_e({\bf p},{\bf q})$. In this paper, first, the detailed structure of ${\cal C}_e({\bf p},{\bf q})$ is well specified and all extreme points are enumerated by a special algorithm. As an application, the exact solution of the minimum-entropy coupling problem is obtained. Second, it is proved that for any strict Schur-concave function $\Psi$ on ${\cal C}(\bf p,q)$, $\Psi$ also takes its minimal value on ${\cal C}_e({\bf p},{\bf q})$. As an application, the exact solution of the minimum-entropy coupling problem is obtained for $(\Phi,\hbar)$-entropy, a large class of entropy including Shannon entropy, R\'enyi entropy and Tsallis entropy etc. Finally, all the above are generalized to multi-marginal case.
\vskip 3mm
{\rm \bf AMS classification (2020):} {\it 94A17, 60E15.}
}


\end{minipage}
\end{center}

 \vskip 5mm
\section{Introduction}
\renewcommand{\theequation}{1.\arabic{equation}}
\setcounter{equation}{0}
The concept of entropy was introduced in thermodynamical and statistical mechanics as a measure of uncertainty or disorganization in a physical system \cite{Bl1,Bl2}. In 1877, L. Boltzmann \cite{Bl2} gave the probabilistic interpretation of
entropy and found the famous formula $S=\kappa \ln W$.
The second law of thermodynamical says that the entropy of a closed system cannot decrease.

\subsection{The minimum-entropy coupling problem}


To reveal the physics of information, C. Shannon \cite{S} introduced the entropy in the communication theory. Let $X$ be a discrete random element with alphabet ${\cal X}$ and probability mass ${\bf p}=\{p(x)=\P(X=x): x\in {\cal X}\}$, the entropy of $X$ (or $\bf p$) is defined by
\be\label{1.1} H(X)=H({\bf p}):=-\sum_{x\in{\cal X}}p(x)\log p(x).\ee

To introduce the minimum-entropy coupling problem, let's first extend the definition of entropy to a pair of random variables.
Let $(X,Y)$ be a two-dimensional random vector in ${\cal X}\times{\cal Y}$ with a joint distribution $P$=$\{p(x,y):x\in{\cal X},y\in{\cal Y}\}$,
the {\it joint entropy} of $(X,Y)$ (or $P$) is defined by
\be\label{1.2}H(X,Y)=H(P):=-\sum_{x\in{\cal X}}\sum_{y\in{\cal Y}}p(x,y)\log p(x,y).\ee An relevant concept in information theory on random vector $(X,Y)$ is the {\it mutual information} (see \cite{CT}, Chapter 2), which is a measure of the amount of information that one random variable contains about the other, and is defined by
\be\label{1.3}I(X,Y):=\sum_{x\in{\cal X}}\sum_{y\in{\cal Y}}p(x,y)\log {\Df {p(x,y)}{p(x)p(y)}},\ee
where $\{p(x):x\in{\cal X}\}$, $\{p(y):y\in{\cal Y}\}$ are the marginal distributions of $X,Y$. By definitions, one has
\be\label{1.4'}I(X,Y)=H(X)+H(Y)-H(X,Y).\ee
Note that in some setting, the {\it maximum} of {\it mutual information} is called the {\it channel capacity}, which plays a key role in information theory through the famous {\it Shannon's second theorem: Channel Coding Theorem} \cite{S}.

For basic concepts and properties in information theory, readers may refer to \cite{CT} and the references therein.


For given marginals $\{p(x):x\in\cal X\}$ and $\{p(y):y\in \cal Y\}$, maximizing $I(X,Y)$ and minimizing $H(X,Y)$ are two sides of a single coin. The problem of finding the minimum-entropy coupling of two discrete probability distribution $\bf p,q$ is called the {\it minimum-entropy coupling problem}.

For integers $m, n\geq 2$, for simplicity, we take ${\cal X}=[m]:=\{1,2,\ldots,m\}$, ${\cal Y}=[n]:=\{1,2,\ldots,n\}$. Note that in the whole paper, $m,n$ always mean integers $\geq 2$. Denote by ${\cal P}_m$, ${\cal P}_n$ the set of all probability distributions on $\cal X$, $\cal Y$ respectively. Clearly, over all ${\bf p}=(p_1,\ldots,p_m)\in{\cal P}_m$, the Shannon entropy of $\bf p$, $H(\bf p)$ takes its minimum $0$ when $\bf p$ is degenerated (i.e. for some $1\leq k\leq m,\ p_k=1$) and takes its maximum $\log m$ when $\bf p$ is uniformly distributed (i.e. $p_k=\frac 1m,\ \forall \ 1\leq k\leq m$). In this sense, entropy is a measure of the uncertainty of a random variable.

For any ${\bf p}\in {\cal P}_m$, ${\bf q}\in {\cal P}_n$, let ${\cal C}(\bf p,q)$ be the set of all couplings (i.e. joint distributions) of ${\bf p,q}$. Clearly ${\cal C}({\bf p},{\bf q})$ forms a $(n-1)\times (m-1)$-dimensional polytope in $\R^{mn}$, denote by ${\cal C}_e({\bf p},{\bf q})$ the vertex set of this polytope, i.e. the set of extreme points of convex set ${\cal C}({\bf p},{\bf q})$. For any $P=\(p_{i,j}\)\in{\cal C}(\bf p,q)$, let's consider its Shannon entropy $H(P)$.
Clearly, in the case when $P$ is the {\it independent} coupling, $H(P)$ takes the maximum $H({\bf p})+H({\bf q})$. The more interesting problem about joint entropy is the following minimum-entropy coupling problem:
\be\label{1} \tilde P:\ \ H(\tilde P)=\inf_{P\in{\cal C}({\bf p},{\bf q})}H(P).\ee The $\tilde P$ which solves the optimization problem (\ref{1}) is called a minimum-entropy coupling.
There should be two main points of concern regarding this issue: the first point is to find out all minimum-entropy couplings and calculate the exact value of the minimal joint entropy; the second point is to specify the structure of a minimum-entropy coupling, a point stems from physicists' interest on the intrinsic ordered structure of systems with minimal entropy. Actually, the minimum-entropy coupling problem (\ref{1}) has already become an important problem in information theory and has been studied deeply in the last two decades, see \cite{CGV,CK,KDVH,KSS,L,MWWC,PRF,R,V,YT} etc.

The natural strategy to solve the minimum-entropy coupling problem can be stated as follows. For any ${\bf p}\in {\cal P}_m$, ${\bf q}\in {\cal P}_n$, by a concave argument, the Shannon entropy H on ${\cal C}(p,q)$ take its minimal value in ${\cal C}_e({\bf p},{\bf q})$, a finite subset of ${\cal C}({\bf p},{\bf q})$. Then the optimization problem (\ref{1}) is transformed to the following optimization problem
\be\label{2} \tilde P:\ \ H(\tilde P)=\min_{P\in{\cal C}_e({\bf p},{\bf q})}H(P).\ee To solve the minimum-entropy coupling problem perfectly, the key is to give a perfect characterization of the extreme point set ${\cal C}_e({\bf p},{\bf q})$. Unfortunately, the structure of ${\cal C}_e({\bf p},{\bf q})$ is complicated enough for general ${\bf p},{\bf q}$, while it is shown in \cite{KSS,V} that this problem is NP-hard, polynomial time approximation algorithms are given in \cite{CGV,CK,KDVH,L,MWW,R} etc.

Recently, depending on the {\it forest} structure of the extreme point, \cite{CK} provided a backtracking algorithm to calculate the minimal joint entropy in exponential time. In the present paper, we shall follow \cite{CK} to finish the complete presentation of the structure of ${\cal C}_e({\bf p},{\bf q})$, then solve the minimum-entropy coupling problem by enumerating all the extreme points with a specified algorithm. In fact, we will introduce a graph representation for the support of a coupling, and then prove that the support of the extreme point possesses a forest structure. Note that our graph is quite different from the graph introduced in \cite{CK}, see Figure~\ref{f2} for an illustration. We emphasis here that, different to the backtracking algorithm provided in \cite{CK}, our algorithm can be successfully generalized to the multi-marginal case.

For the special structure of a minimum-entropy coupling $\tilde P$, it is shown recently in \cite{MWW1} that, for any ${\bf p},{\bf q}\in{\cal P}_n$, $\tilde P$ is {\it essentially order-preserving}. Note that this in some sense fulfills the gap for us to interpret entropy as a measure of {\it system disorder}. In the present paper, besides the {\it forest } structure, more special structures of a minimum-entropy coupling are revealed (see Theorem~\ref{t2} and Figure~\ref{f1}).


%
%
%
%
%

Note that inferring an unknown joint distribution of two random
variables with given marginals is an old problem in the area of probabilistic inference. As far as we know, the problem may go back at least to Frechet \cite{F} and Hoeffding \cite{H}, who
studied the question of identifying the extremal joint distribution that maximizes (resp., minimizes) their correlation, for more literatures in this area and more applications in pure and applied sciences, readers may refer to \cite{BS,CFR,DKS,LDKH} etc.

\subsection{Statement of the result}

Recall that a permutation $\sigma$ is a bijective map from $[m]$ into itself, denote by $\Sigma_m$ the set of all permutations. For any ${\bf p}=(p_1,p_2,\ldots,p_m)\in{\cal P}_m$, define $\sigma{\bf p}:=(p_{\sigma(1)},p_{\sigma(2)},\ldots,p_{\sigma(m)})$ and denote by $\bar{\bf p}$ the permutation of $\bf p$ such that $\bar p_1\geq\bar p_2\geq\ldots\geq \bar p_m$. By the definition (\ref{1.1}), one has
\be\label{1.3}H({\bf p})=H(\sigma{\bf p}),\ \forall\ \sigma\in \Sigma_m,\ee i.e. the Shannon entropy is a symmetric function.
For random variable $X$ with distribution $\bf p$, random variable $\sigma X$ has the distribution $\sigma^{-1} \bf p$, where $\sigma^{-1}$ is the inverse of $\sigma$.

For each ${\bf p}\in{\cal P}_m$, let $F_{\bf p}$ be the cumulative distribution function defined by
\be\label{1.4}F_{\bf p}(i):=\sum_{k=1}^ip_k,\ 1\leq i\leq n.\ee

For any ${\bf p}\in {\cal P}_m$, ${\bf q}\in {\cal P}_n$, $P\in{\cal C}(\bf p,q)$, suppose random vector $(X,Y)$ is distributed according to $P$. For any permutation pair $(\sigma,\pi)\in \Sigma_m\times\Sigma_n$, denote by $P(\sigma,\pi)$ the joint distribution of $(\sigma X,\pi Y)$, then
\be\label{1.5} P(\sigma,\pi)\in{\cal C}({\sigma^{-1}{\bf p},\pi^{-1}{\bf q}})\ {\rm and}\ \ H(P(\sigma,\pi))=H(P).\ee

For any $m\geq 2$, let ${\cal P}^+_m=\{{\bf p}\in{\cal P}_m: p_k>0,\forall \ 1\leq k\leq m\}$. In this paper, we shall study the optimization problem (\ref{1}) for ${\bf p}\in{\cal P}^+_m$, ${\bf q}\in{\cal P}^+_n$.

\bd\label{d0}  For any ${\bf p}\in {\cal P}^+_m$, ${\bf q}\in {\cal P}^+_n$, define the \textbf{structure constant} of pair $(\bf p,q)$ as the following.
\be\label{1.6}\kappa({\bf p,q}):=\max_{(\sigma,\pi)\in \Sigma_m\times\Sigma_n}|\{F_{\sigma \bf p}(i):1\leq i<m\}\cap\{F_{\pi \bf q}(j):1\leq j<n\}|+1.
\ee
\ed
\br We call $\kappa(\bf p,q)$ the {\it structure constant}, since for any $1\leq k\leq \kappa(\bf p,q)$, there exists $P\in {\cal C}(\bf p,q)$ such that $P(\sigma,\pi)$ possesses the block structure as given in (\ref{2.6}) for some permutation pair $(\sigma,\pi)$. \er

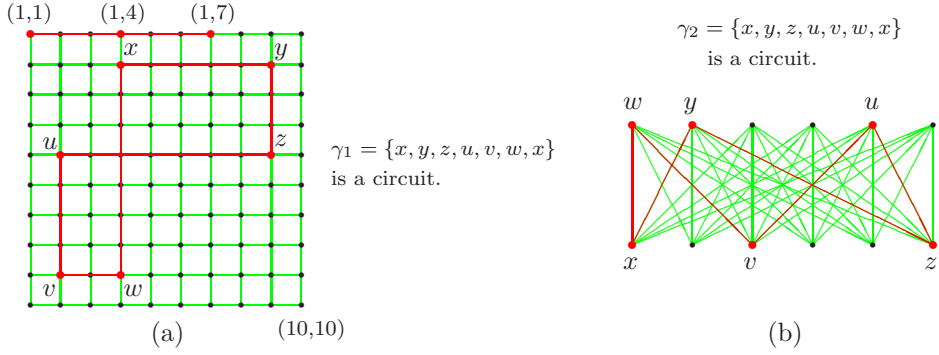
\begin{figure}
\unitlength=0.4mm
\begin{picture}(100,100)(-50,0)
\put(0,0){\color{green}\line(1,0){90}}\put(0,0){\color{green}\line(0,1){90}}
\put(0,10){\color{green}\line(1,0){90}}\put(10,0){\color{green}\line(0,1){90}}
\put(0,20){\color{green}\line(1,0){90}}\put(20,0){\color{green}\line(0,1){90}}
\put(0,30){\color{green}\line(1,0){90}}\put(30,0){\color{green}\line(0,1){90}}
\put(0,40){\color{green}\line(1,0){90}}\put(40,0){\color{green}\line(0,1){90}}
\put(0,50){\color{green}\line(1,0){90}}\put(50,0){\color{green}\line(0,1){90}}
\put(0,60){\color{green}\line(1,0){90}}\put(60,0){\color{green}\line(0,1){90}}
\put(0,70){\color{green}\line(1,0){90}}\put(70,0){\color{green}\line(0,1){90}}
\put(0,80){\color{green}\line(1,0){90}}\put(80,0){\color{green}\line(0,1){90}}
\put(0,90){\color{green}\line(1,0){90}}\put(90,0){\color{green}\line(0,1){90}}
\multiput(0,0)(10,0){10}{\circle*{2}}\multiput(0,20)(10,0){10}{\circle*{2}}
\multiput(0,10)(10,0){10}{\circle*{2}}\multiput(0,30)(10,0){10}{\circle*{2}}
\multiput(0,40)(10,0){10}{\circle*{2}}\multiput(0,50)(10,0){10}{\circle*{2}}
\multiput(0,60)(10,0){10}{\circle*{2}}\multiput(0,70)(10,0){10}{\circle*{2}}
\multiput(0,80)(10,0){10}{\circle*{2}}\multiput(0,90)(10,0){10}{\circle*{2}}
\put(0,90){\color{red}\circle*{3}}\put(30,90){\color{red}\circle*{3}}\put(60,90){\color{red}\circle*{3}}
\put(-8,95){\footnotesize (1,1)}\put(82,-10){\footnotesize (10,10)}\put(23,95){\footnotesize (1,4)}\put(53,95){\footnotesize (1,7)}
\put(0,90){\color{red}\line(1,0){30}}\put(30,90){\color{red}\line(1,0){30}}
\put(30,80){\color{red}\circle*{3}}\put(80,80){\color{red}\circle*{3}}\put(80,50){\color{red}\circle*{3}}\put(10,50){\color{red}\circle*{3}}
\put(10,10){\color{red}\circle*{3}}\put(30,10){\color{red}\circle*{3}}
\put(10,10){\color{red}\line(1,0){20}}\put(10,50){\color{red}\line(1,0){70}}\put(10,10){\color{red}\line(1,0){20}}\put(30,80){\color{red}\line(1,0){50}}
\put(10,10){\color{red}\line(0,1){40}}\put(30,10){\color{red}\line(0,1){70}}\put(80,50){\color{red}\line(0,1){30}}
\put(4,52){$u$}\put(4,3){$v$}\put(31,3){$w$}\put(31,83){$x$}\put(81,83){$y$}\put(81,53){$z$}
\put(100,50){\footnotesize $\gamma_1=\{x,y,z,u,v,w,x\}$}
\put(100,40){\footnotesize {\rm is a circuit}.}
\put(40,-12){(a)}\put(245,-12){(b)}

\multiput(200,20)(20,0){6}{\color{green}\line(0,1){40}}
\multiput(200,20)(20,0){5}{\color{green}\line(1,2){20}}\multiput(200,20)(20,0){4}{\color{green}\line(1,1){40}}
\multiput(200,20)(20,0){3}{\color{green}\line(3,2){60}}
\multiput(200,20)(20,0){2}{\color{green}\line(2,1){80}}
\multiput(200,20)(20,0){1}{\color{green}\line(5,2){100}}
\multiput(200,60)(20,0){5}{\color{green}\line(1,-2){20}}\multiput(200,60)(20,0){4}{\color{green}\line(1,-1){40}}
\multiput(200,60)(20,0){3}{\color{green}\line(3,-2){60}}
\multiput(200,60)(20,0){2}{\color{green}\line(2,-1){80}}
\multiput(200,60)(20,0){1}{\color{green}\line(5,-2){100}}
\multiput(200,20)(20,0){6}{\circle*{2}}\multiput(200,60)(20,0){6}{\circle*{2}}
\put(200,20){\color{red}\line(1,2){20}}\put(220,60){\color{red}\line(2,-1){80}}
\put(280,60){\color{red}\line(1,-2){20}}\put(240,20){\color{red}\line(1,1){40}}
\put(200,60){\color{red}\line(1,-1){40}}\put(200,20){\color{red}\line(0,1){40}}
\put(197,12){$x$}\put(217,65){$y$}\put(297,12){$z$}\put(277,65){$u$}\put(237,12){$v$}\put(197,65){$w$}
\put(200,20){\color{red}\circle*{3}}\put(240,20){\color{red}\circle*{3}}\put(300,20){\color{red}\circle*{3}}
\put(200,60){\color{red}\circle*{3}}\put(220,60){\color{red}\circle*{3}}\put(280,60){\color{red}\circle*{3}}
\put(215,90){\footnotesize $\gamma_2=\{x,y,z,u,v,w,x\}$}
\put(215,80){\footnotesize {\rm \ \ \ \ is a circuit}.}

\end{picture}
\vskip 0mm
\begin{center}
\begin{minipage}{13.5cm}

\caption{\it\small (a) is the graph $G_{m,n}$ with $m=n=10$, it has $n^2$ vertices and $2n\binom{n}{2}=n^3-n$ edges. In $G_{10,10}$, $\gamma_1$ is a circuit, but $\gamma=\{(1,1),(1,4),(1,7),(1,1)\}$ is not a circuit. (b) is the graph $\bar G_{m,n}$ introduced by \cite{CK} with $m=n=6$, which is a bipartite graph with $2n$ vertices and $n^2$ edges, where $\gamma_2$ is a circuit in $\bar G_{m,n}$.
}\label{f2}
\end{minipage}
\end{center}

\end{figure}

\bd \label{d1} For any $m,n\geq 2$, let $G_{m,n}=(V_{m,n},E_{m,n})$ be the graph with vertex set $V_{m,n}=[m]\times[n]$ and edge set $$E_{m,n}:=\{\langle u,v\rangle:u=(u_1,u_2),v=(v_1,v_2)\in V_{m,n}, u\not= v {\rm \ and \ } |u_1-v_1|\cdot|u_2-v_2|=0\}.$$
As a basic concept in graph theory \cite{B}, a sequence $\gamma=\{(i_k,j_k):0\leq k\leq s\}$ of points in $G_{m,n}$
is called a \textbf{path}, if $\langle (i_k,j_k),(i_{k+1},j_{k+1})\rangle$ $\in E_{m,n}$ for all $0\leq k\leq s-1$.
Particularly, in the present paper, a path $\gamma$ is called \textbf{directed}, if $i_k\leq i_{k+1},j_k\leq j_{k+1}$ for all $0\leq k\leq s-1$. A path $\gamma$ is called \textbf{continues}, if $|i_{k+1}-i_{k}|+|j_{k+1}-j_k|=1$ for all $0\leq k\leq s-1$.
A path $\gamma$ is called a \textbf{circuit}, if $(i_0,j_0)=(i_s,j_s)$, $(i_k,j_k)\neq(i_l,j_l)$ for all $0\leq k<l\leq s-1$ and
$$\Dp\prod_{k=0}^{s}|(i_{k+2}-i_k)(j_{k+2}-j_k)|>0,\ {\rm with }\ (i_{s+t},j_{s+t})=(i_s,j_s),\  t=1,2.$$
For any $V\subset V_{m,n}$, see $V$ as the subgraph of $G_{m,n}$ with vertex set $V$ and edge set $E_V=\{\langle u,v\rangle\in E_{m,n}:u,v\in V\}$. A path $\gamma$ is called a path in $V$, if each vertex of $\gamma$ lies in $V$. $V$ is called a \textbf{forest}, if there is no circuit in $V$. A forest $V$ is called a \textbf{tree}, if it is connected, i.e., for any distinct $(i,j),(i',j')\in V$, there exists a path $\gamma=\{(i_k,j_k):0\leq k\leq s\}$ in $V$ such that $(i_0,j_0)=(i,j),\ (i_s,j_s)=(i',j')$. $V$ is called \textbf{complete}, if $\{i:(i,j)\in V\}=[m]$, $\{j:(i,j)\in V\}=[n]$.
\ed
\br The circuit defined above is not the same as what defined in classic graph theory, see \cite{B}. For an example, $\gamma=\{(1,1),(1,4),(1,7),(1,1)\}$ is a circuit in classic significance, but it is not a circuit according to the above Definition~\ref{d1}, see Figure~\ref{f2}.\er

For any $V\subset V_{m,n}$, let $A=A(V)=(a_{i,j})_{m\times n}$ be the {\it indicator matrix} of $V$ such that $a_{i,j}=I_{(i,j)\in V}$, where $I_{(i,j)\in V}$ is the indicator function, namely
$$I_{(i,j)\in V}=\left\{\ba{ll}1,&\ {\rm if}\ (i,j)\in V;\\[2mm] 0,&\ {\rm otherwise.}\ea\right.$$ Let $\bar P=P(V)=\frac 1{|V|}A(V)$ be a probability matrix. Using Lemma~\ref{l3} and Theorem~\ref{t3} to $\bar P$, we have

\bp\label{p0}Suppose $V\subset V_{m,n}$, then
\begin{description}
  \item[\ \ \ \ \it i)] if $V$ is a forest with $k$ connected components, then $V$ is complete if and only if $|V|=m+n-k$;
  \item[\ \ \ \ \it ii)] if $|V|=m+n-1$, then $V$ is a tree if and only if $V$ is complete and connected;
  \item[\ \ \ \ \it iii)] if $V$ is complete and connected, then $|V|\geq m+n-1$ and $|V|=m+n-1$ if and only if $V$ is a tree.
\end{description}

\ep

For any nonnegative matrix $A=(a_{i,j})_{m\times n}$ (i.e. all its entries are nonnegative), let $V(A)=\{(i,j):a_{i,j}\not=0\}$ be the support of $A$. Write $V(A)$ as the disjoint union of the following $V_s(A), s=1,2,3$: $$V_1(A)=\left\{(i,j)\in V(A):\sum_{k=1}^m a_{k,j}=\sum_{l=1}^n a_{i,l}=a_{i,j}\right\},$$
$$V_2(A)=\left\{(i,j)\notin V_1(A):\sum_{k=1}^ma_{k,j}\ {\rm or} \ \sum_{l=1}^n a_{i,l}=a_{i,j}\right\}$$ and $V_3(A)=V(A)\setminus \(V_1(A)\cup V_2(A)\)$. Let $$V^r_2(A):=\left\{(i,j)\in V_2(A):\sum^n_{l=1}a_{i,l}=a_{i,j}\right\}$$ and $V^c_2(A):=V_2(A)\setminus V^r_2(A)$. Note that from the view of a passenger walking alone a path in graph $V(A)$, $V_1(A)$ is the {\it isolated vertex} set, $V_2(A)$ is the {\it row or column passable vertex} set and $V_3(A)$ is the {\it turning vertex} set of $A$. As a basic fact, we declare the following proposition without proof.

\bp\label{p-} For any ${\bf p}\in{\cal P}^+_m$, ${\bf q}\in {\cal P}^+_n$, for any $P\in {\cal C}({\bf p,q})$, $V(P)$ is complete; if furthermore $\kappa({\bf p,q})=1$, then $V(P)$ is complete and connected.\ep

\br The graph introduced in \cite{CK} is $\bar G_{m,n}:=(V,E)$, where $V=V_{r}\cup V_{c}$ with $|V_{r}|=m$, $|V_{c}|=n$, $E=\{\langle i,j\rangle:i\in V_{r}, j\in V_{c}\},$ see Figure~\ref{f2} (b). For any probability matrix $P=(p_{i,j})_{m\times n}$, while we define the  subgraph $V(P)$ of $G_{m,n}$, S. Compton etc. \cite{CK} have defined the subgraph $E(P)=(V,E(P))$ of $\bar G_{m,n}$ with $V:=V_{r}(P)\cup V_{c}(P)$$=\{i\in V_{r}:$ for some $j\in V_{c}$, $p_{i,j}>0\}\cup$$\{j\in V_{c}:$ for some $i\in V_{r}$, $p_{i,j}>0\}$ and $E(P)=\{\langle i,j\rangle:p_{i,j}>0\}$. $V(P)$ and $E(P)$ are quite different, and obviously $V(P)$ possesses a more detailed structure. $V(P)$ and $E(P)$ are associated by the following property: there exists a circuit in $V(P)$ if and only if there exists a circuit in $E(P)$. \er

\bl\label{l0} For any nonnegative matrix $A=(a_{i,j})_{m\times n}$, $m,n\geq 2$, if $V(A)\subset V_{m,n}$ is a forest, then $V_1(A)\cup V_2(A)\not=\emptyset$, $V_3(A)$ is a forest whenever $V_3(A)\not=\emptyset$; if furthermore $V(A)$ is a tree, then $V_1(A)=\emptyset$, $V_2(A)\not=\emptyset$ and $V_3(A)$ is a tree whenever $V_3(A)\not=\emptyset$.\el

{\it Proof.} If $V(A)=V_3(A)$, for any $(i_0,j_0)\in V_3(A)=V(A)$, for any $s\geq 2$, we can find a path $\gamma=\{(i_k,j_k)\in V_3(A):0\leq k\leq s\}$ such that $|i_{k+2}-i_k|\cdot|j_{k+2}-j_k|>0, \ {\rm for\ all}\ 0\leq k\leq s-2$. Since $V(A)$ is a forest, i.e. there is no circuit in $V(A)$, then all vertexes in $\gamma$ are distinct, this implies that $|V_{3}(A)|=|V(A)|\geq s$, a contradiction to the arbitrariness of $s$.

In the case when $V(A)$ is a tree, by definition, $V_1(A)=\emptyset$. If $V_3(A)$ is not a tree, then there exists $(i,j), (i',j')\in V_3(A)$, $(i'',j'')\in V_2(A)$ such that $\{(i,j), (i',j')\}$ does not form a path, but $\{(i,j), (i'',j'')$, $(i',j')\}$ forms a path. By Definition~\ref{d1}, this implies that $(i'',j'')\in V_3(A)$, a contradiction.\QED

Denote ${\cal T}=\{T\subset V_{m,n}:$ $T$ is a tree and $|T|=m+n-1\}$. By Proposition~\ref{p0}, for any $T\in {\cal T}$, $T$ is complete.

\bd\label{1.3}For any ${\bf p}\in {\cal P}^+_m$, ${\bf q}\in {\cal P}^+_n$, suppose $T\in{\cal T}$ and $P\in{\cal C}(\bf p,q )$. $P$ is called \textbf{consistent} with $T$, if $V(P)\subset T$. \ed

The following proposition will play a key role in the description of the extreme points set ${\cal C}_e({\bf p},{\bf q})$.
\bp\label{p1}For any ${\bf p}\in {\cal P}^+_m$, ${\bf q}\in {\cal P}^+_n$ and for any $T\in{\cal T}$, there exists at most one $P\in{\cal C}(\bf p,q)$ such that $P$ is consistent with $T$.\ep

{\it Proof.} It suffices to prove that: for any $T\in{\cal T}$, there exists at most one $P\in{\cal C}(\bf p,q)$ such that $p_{i,j}=0$ for all $(i,j)\notin T$.

For any tree $T\in{\cal T}$, let $A=A(T)=(a_{i,j})_{m\times n}$ be the indicator matrix of $T$. By Lemma~\ref{l0}, $V_1(A)=\emptyset$, $V_2(A)\not=\emptyset$ and $T=V(A)=V_2(A)\cup V_3(A).$ Recall that $V_2(A)=V^r_2(A)\cup V^c_2(A)$.

Now, if $P$ is a probability matrix in ${\cal C}(\bf p,q)$ such that $p_{i,j}=0$ for all $(i,j)\notin T$, then for $(i,j)\in V_2(A)$,
\be\label{1.11}p_{i,j}=\left\{\ba{ll}p_i,&{\rm if\ } (i,j)\in V^r_2(A);\\[8mm]
q_j,&{\rm if\ }(i,j)\in V^c_2(A).\ea\right.\ee

If $V_3(A)=\emptyset$, then we finish the definition of $P$. Otherwise, Let $A_1$ be the submatrix of $A$ such that $V(A_1)=V_3(A). $ Let $P_1$ be the submatrix of $P$ satisfies: $p_{i,j}$ is an entry in $P_1$ if and only if $a_{i,j}$ is an entry in $A_1$. For any entry $a_{i,j}$ of $A_1$, let
\be\label{1.12}p^1_i=p_i-\sum_{l:(i,l)\in V^c_2(A)}p_{i,l}, \ \ q^1_j=q_j-\sum_{k:(k,j)\in V^r_2(A)}p_{k,j}\ee be the corresponding row and column summations of $P_1$.

By Lemma~\ref{l0},
$V(A_1)$ is a tree and $V(A_1)=V_2(A_1)\cup V_3(A_1)$, $V_2(A_1)=V^r_2(A_1)\cup V^c_2(A_1)\not=\emptyset.$ Then, for any $(i,j)\in V_2(A_1)$,

\be\label{1.13}p_{i,j}=\left\{\ba{ll}p^1_i,&{\rm if\ } (i,j)\in A^r_2(A_1);\\[8mm]
q^1_j,&{\rm if\ }(i,j)\in A^c_2(A_1).\ea\right.\ee

Repeat the above procedure for $\xi\geq 2$ until $V_3(A_\xi)=\emptyset$: If $V_3(A_{\xi-1})\not=\emptyset$, let $A_\xi$ be the submatrix of $A_{\xi-1}$ such that $V(A_\xi)=V_3(A_{\xi-1}). $ Let $P_\xi$ be the submatrix of $P_{\xi-1}$ satisfies: $p_{i,j}$ is an entry in $P_\xi$ if and only if $a_{i,j}$ is an entry in $A_\xi$. For any entry $a_{i,j}$ of $A_\xi$, let
\be\label{1.7}p^\xi_i=p^{\xi-1}_i-\sum_{l:(i,l)\in V^c_2(A_{\xi-1})}p_{i,l},\ \ q^\xi_j=q^{\xi-1}_j-\sum_{k:(k,j)\in V^r_2(A)}p_{k,j}.\ee

By Lemma~\ref{l0},
$V(A_\xi)$ is a tree and $V(A_\xi)=V_2(A_\xi)\cup V_3(A_\xi)$, $V_2(A_\xi)=V^r_2(A_\xi)\cup V^c_2(A_\xi)\not=\emptyset.$ Then, for any $(i,j)\in V_2(A_{\xi})$,

\be\label{1.8}p_{i,j}=\left\{\ba{ll}p^\xi_i,&{\rm if\ } (i,j)\in A^r_2(A_\xi);\\[8mm]
q^\xi_j,&{\rm if\ }(i,j)\in A^c_2(A_\xi).\ea\right. \ee

Let $\xi_0:=\min\{\xi\geq 0:V_3(A_\xi)=\emptyset\}$, then
$$T=\bigcup_{s=0}^{\xi_0}V_2(A_s)\ {\rm (with}\ A_0=A{\rm )}.$$ Thus
$p_{i,j}$ is determined by (\ref{1.11}), (\ref{1.13}) and (\ref{1.8}) for all $(i,j)\in T$. Namely, $P$ is uniquely determined by $T$ and $\bf p,q$. Obviously, $P$ is consistent with $T$.\QED

\br\label{r1} The above proof actually provides an algorithm to obtain the unique $P$ whenever it exists. Note that $P$ exists if and any if $p_{i,j}$ defined in the proof of Proposition~\ref{p1} is always nonnegative. Actually, by Lemma~\ref{l4}, we can obtain the unique $P$ by solving a system of linear equations. \er

\br For ${\bf p,q}$ with $\kappa({\bf p,q})=1$, by Proposition~\ref{p0}, iii), Propositions~\ref{p-} and \ref{p1}, for any $V\subset[m]\times[n]$ with $|V|=m+n-1$, there exists at most one $P\in{\cal C}(\bf p,q)$ such that $P$ is consistent with $V$.\er
Let ${\hC}({\bf p,q})=\{P\in {\cal C}(\bf p,q):$ there exists $T\in {\cal T}$ such that $P$ is consistent with $T\}$. By Proposition~\ref{p1}, $|\hC({\bf p,q})|\leq |{\cal T}|\leq \binom {mn}{m+n-1}<\infty$. In the whole paper we write $C(m,n):=\binom {mn}{m+n-1}$.

We introduce our {\bf Main Theorem} as the following.
\bt\label{t1} For any ${\bf p}\in {\cal P}^+_m$, ${\bf q}\in {\cal P}^+_n$, one has
\be\label{1.9'}{\cal C}_e({\bf p},{\bf q})={\hC}({\bf p},{\bf q}).\ee Thus, if $\tilde P\in{\cal C}(\bf p,q)$ solves the optimization problem (\ref{1}), then $\tilde P\in{\mathscr C}(\bf p,q)$ and
\be\label{1.9} H(\tilde P)=\min_{P\in \hC({\bf p},{\bf q})}H(P).\ee \et

As a corollary of the main Theorem~\ref{t1}, the minimum-entropy coupling problem for R\'enyi entropy \cite{Re} and Tsallis entropy \cite{Ts} can be similarly addressed. Note that for parameter $\a$, $\a\geq 0,\ \a\not=1$, the R\'enyi entropy and the Tsalls entropy are defined by
\be\label{4.5}H({\bf p})=H^R_\a({\bf p}):=\frac 1{1-\a}\log \(\sum_{i=1}^m p_i^\a\),\ {\bf p}\in{\cal P}_m \ee
and
\be\label{4.6} H({\bf p})=H^T_\a({\bf p}):=\frac 1{1-\a}\(\sum_{i=1}^m p_i^\a-1\),\ {\bf p}\in{\cal P}_m \ee respectively. It is straightforward to check that the R\'enyi entropy and the Tsallis entropy are all strictly concave functions on ${\cal C}({\bf p},{\bf q})$.

\bco\label{cor1} For any ${\bf p}\in {\cal P}^+_m$, ${\bf q}\in {\cal P}^+_n$, if $\tilde P\in{\cal C}(\bf p,q)$ solves the optimization problem (\ref{1}) for R\'enyi entropy or Tsallis entropy, then $\tilde P\in{\mathscr C}({\bf p},{\bf q})={\cal C}_e({\bf p},{\bf q})$ and
\be\label{1.10} H(\tilde P)=\min_{P\in \hC({\bf p},{\bf q})}H(P).\ee \eco

\vskip 5mm
\section{Structure of the minimum-entropy couplings and proof of the Main Theorem}
\renewcommand{\theequation}{2.\arabic{equation}}
\setcounter{equation}{0}
Suppose $A=(a_{i,j})_{m\times n}$ is a nonnegative matrix such that $C:=\sum_{i=1}^m\sum_{j=1}^n a_{i,j}>0$. We generalize the definition of entropy for nonnegative matrix $A$ as
\be\label{2.1}H(A):=-\sum_{i=1}^m\sum_{j=1}^n a_{i,j}\log a_{i,j}.\ee Let $P=C^{-1}A$, a probability matrix, then
\be\label{2.2} H(A)=C H(P)-C\log C.\ee

In this section, we will use the following local optimization lemmas developed in \cite{MWWC} to study the special structure of a minimum-entropy coupling.

\noindent{\bf Lemma 1}[Lemma 2.2 in \cite{MWWC}]. {\it For any second order nonnegative matrix $A=\(a_{i,j}\)_{2\times 2}$. Suppose that $a_{1,1}\vee a_{2,2}\geq a_{1,2}\vee a_{2,1}$, denote $b=a_{1,2}\wedge a_{2,1}$. Let $A'=\(a'_{i,j}\)_{2\times 2}$ such that $a'_{i,i}=a_{i,i}+b,\ i=1,2$, $a'_{i,j}=a_{i,j}-b,\ i\neq j.$ Then $H(A)\geq H(A')$. Furthermore, if $b>0$, then $H(A)>H(A')$. Where $\cdot\vee\cdot$, $\cdot\wedge\cdot$ means $\max\{\cdot,\cdot\}$, $\min\{\cdot,\cdot\}$ respectively.}

\noindent{\bf Lemma 2}{[\rm Lemma\ 2.3\ in \cite{MWWC}]}. {\it For any second order nonnegative matrix $A=\(a_{i,j}\)_{2\times 2}$. Suppose that $a_{1,1}+a_{1,2}\geq a_{2,1}+a_{2,2}, a_{1,1}+a_{2,1}\geq a_{1,2}+a_{2,2}$ and $a_{1,1}+a_{1,2}\geq a_{1,1}+a_{2,1}$. Let $b=a_{1,2}\wedge a_{2,1}$, define $A'$ as in Lemma~1, then $H(A)\geq H(A')$.}

As a consequence of Lemmas~1 and 2, we introduce an additional lemma for local optimization as the following.

\bl\label{l1} For any $2\times n$ nonnegative matrix $A=(a_{i,j})$ with $a_{2,k}=0$, $2\leq k\leq n$, let
$A'$=$(a'_{i,j})$ be the $2\times n$ matrix such that $a'_{1,1}=\sum_{k=1}^n a_{1,k}$, $a'_{1,k}=0,\ 2\leq k\leq n$; $a'_{2,1}=a_{2,1}-\sum_{k=2}^n a_{1,k}$, $a'_{2,k}=a_{1,k},\ 2\leq k\leq n$. i.e.
$$ A=\(\ba{llll}\Dp a_{1,1}&a_{1,2}&\ldots&a_{1,n}\\[3mm]
a_{2,1}&0&\ldots&0\ea\),\ \ A'=\(\ba{llll}\sum_{k=1}^n a_{1,k}&0&\ldots&0\\[3mm]
a_{2,1}-\sum_{k=2}^n a_{1,k}&a_{1,2}&\ldots&a_{1,n}\ea\).$$ If $\sum_{k=2}^na_{1,k}\leq a_{2,1}\leq \sum_{k=1}^n a_{1,k}$, then $H(A)\geq H(A')$.
\el

By Lemmas~1, 2 and Lemma~\ref{l1}, we obtain the following local optimization theorem.

\noindent{\bf Theorem 3}[Theorem 2.5 in \cite{MWWC}]. {\it Suppose ${\bf p}\in {\cal P}^+_m$, ${\bf q}\in {\cal P}^+_n$. Let
$A$ be the submatrix of $P$ which satisfies the conditions in Lemma 1, Lemma 2 or Lemma~\ref{l1}, and $A'$ be the corresponding matrix of $A$. Let $P'$ be the matrix obtained from $P$ by transforming $A$ to $A'$, then $P'\in {\cal C}(\bf p,q)$ and $H(P)\geq H(P')$. In particular, $H(P)>H(P')$ if and only if $H(A)>H(A')$.}
\bd\label{d2} For any ${\bf p}\in {\cal P}^+_m$, ${\bf q}\in {\cal P}^+_n$, $P\in {\cal C}(\bf p,q)$ is called \textbf{local optimal}, if it can not be further optimised by Lemma~1, Lemma~2 and Lemma~\ref{l1}.
If $\tilde P\in{\cal C}(\bf p,q)$ solves the optimization problem (\ref{1}), i.e. $\tilde P$  is a minimum-entropy coupling, then $\tilde P$ is {\it local optimal}.\ed
As the main result in \cite{MWW1}, for $m=n$, it is proved that, if $\tilde P$ is a minimum-entropy coupling and random variable $(X,Y)$ is distributed according to $\tilde P$, then there exists permutation pair $(\sigma,\pi)\in\Sigma_m\times\Sigma_m$ such that
\be\label{2.11}\P(\sigma X\leq\pi Y)=1.\ee In this sense, $(X,Y)$ or $\tilde P$ is called {\it essentially order-preserving}. Note that equation (\ref{2.11}) is equivalent to the {\it upper triangular} structure of $\tilde P(\sigma,\pi)$, the distribution of $(\sigma X,\pi Y)$, and then equivalent to the fact $F_{\pi^{-1}{\bf q}}\leq F_{\sigma^{-1}{\bf p}}$, i.e. $\pi^{-1}{\bf q}$ is majorized by $\sigma^{-1}{\bf p}$.

In this section, we try to reveal more detailed structures of a local optimal coupling, see the following Theorem~\ref{t2} and Theorem~\ref{t3}, and these structures will play key roles in the proof of Theorem~\ref{t1}.
\bt\label{t2}For any ${\bf p}\in {\cal P}^+_m$, ${\bf q}\in {\cal P}^+_n$, $P\in {\cal C}(\bf p,q)$ is local optimal if and only if the following hold
\begin{description}
  \item[\ \ \ \ \ \ 1.] $V(P)$ is a complete forest; furthermore, if additionally $\kappa(\bf p,q)=1$, then $V(P)$ is a complete tree.
  \item[\ \ \ \ \ \ 2.] For any $2\times n_1$ submatrix $A$ of $P$, $2\leq n_1\leq n$ such that all entries in one row are positive, and only one entry in the other row is positive, without loss of generality, suppose
  $$A=\(\ba{llll}a_{1,1}&a_{1,2}&\ldots&a_{1,n_1}\\[3mm]
  a_{2,1}&0&\ldots &0 \ea\).$$ Then either
  \begin{itemize}
    \item $a_{2,1}\geq \Dp\sum_{k=1}^{n_1}a_{1,k}$, or
    \item $a_{1,1}=\min\{a_{1,k}:1\leq k\leq n_1\}$ and $a_{1,k}\geq a_{1,1}+a_{2,1}$ for all $2\leq k\leq n_1$.
  \end{itemize}
\item[\ \ \ \ \ \ 3.] The above item 2. holds for $P^{T}$, the transpose of $P$.
\end{description}
\et
{\it Proof.} By Definition~\ref{d2}, it is only necessary to prove that, for any local optimal $P\in{\cal C}(\bf p,q)$, $V(P)$ is a forest. Although \cite{CK}
has given a proof for the case of Shannon entropy, we still give a proof based on the previous lemmas. We point out that our proof can be successfully extended to the general Schur-concave function case.

First of all, by Lemma~1, for any $2$-nd order submatrix $A$ of $P$, at least one entry of $A$ is zero. Now, if $\gamma=\{(i_k,j_k):0\leq k\leq s\}$ is a circuit in $V(P)$, suppose that $p_{i_{k_0},j_{k_0}}=\max\{p_{i_k,j_k}:0\leq k\leq s\}$. Without  loss  of generality, assume that $0<k_0<s-1$ and $i_{k_0+1}=i_{k_0}$, $j_{k_0+1}>j_{k_0}$, $i_{k_0}<i_{k_0-1}$, $j_{k_0}=j_{k_0-1}$. Let's consider the following $2$-nd order submatrix of $P$:
 $$A=\(\ba{ll}p_{i_{k_0},j_{k_0}}& p_{i_{k_0+1},j_{k_0+1}}\\[5mm]
p_{i_{k_0-1},j_{k_0-1}}&p_{i_{k_0-1},j_{k_0+1}}\ea\).$$
By the argument mentioned above and the definition of a circuit, one has $p_{i_{k_0-1},j_{k_0+1}}=0$, $p_{i_{k_0},j_{k_0}}\geq p_{i_{k_0+1},j_{k_0+1}}\vee p_{i_{k_0-1},j_{k_0-1}}\geq b:= p_{i_{k_0+1},j_{k_0+1}}\wedge p_{i_{k_0-1},j_{k_0-1}}>0$. Let
$$A'=\(\ba{ll}p_{i_{k_0},j_{k_0}}+b& p_{i_{k_0+1},j_{k_0+1}}-b\\[5mm]
p_{i_{k_0-1},j_{k_0-1}}-b&\hskip 8mm b\ea\),$$
then by Lemma~1, $H(A)>H(A')$. Let $P'\in {\cal C}(\bf p,q)$ be the probability matrix obtained from $P$ by $A'$ taking the place of $A$, by Theorem~3, one has $H(P)>H(P')$, a contradiction to Definition~\ref{d2}. So, there is no circuit in $V(P)$ and $V(P)$ is a forest.

Finally, if $V(P)$ is a forest but not a tree, then there exists some permutation pair $(\sigma,\pi)\in \Sigma_m\times\Sigma_n$ such that $P(\sigma,\pi)\in{\cal C}(\sigma^{-1}{\bf p},\pi^{-1}{\bf q})$ has the following block structure
$$P(\sigma,\pi)=\(\ba{ll}P_{1}&\ 0\\[3mm]
\,0&\ P_{2}\ea\),$$ where $P_{1}$ (resp. $P_{2}$) is a $m_1\times n_1$ (resp. $(m-m_1)\times(n-n_1)$) nonnegative matrix for some $1\leq m_1< m,\ 1\leq n_1<n$. The $0$'s are the corresponding zero matrixes. This implies that
$$\{F_{\sigma^{-1} \bf p}(i):1\leq i<m\}\cap\{F_{\pi^{-1} \bf q}(j):1\leq j<n\}\not=\emptyset$$ and then $\kappa(\bf p,q)>1$, a contradiction.
\QED


\bt\label{t3}For any ${\bf p}\in {\cal P}^+_m$, ${\bf q}\in {\cal P}^+_n$, if $P\in {\cal C}(\bf p,q)$ is local optimal, then
\be\label{2.3} m+n-\kappa({\bf p,q})\leq |V(P)|\leq m+n-1.\ee
\et

Before giving a proof to Theorem~\ref{t3}, we introduce the following lemma.

%
%
%

\bl\label{l3}For any ${\bf p}\in {\cal P}^+_m$, ${\bf q}\in {\cal P}^+_n$, if $P\in {\cal C}(\bf p,q)$ is local optimal and $V(P)$ is a tree, then
\be\label{2.5} |V(P)|= m+n-1.\ee \el

{\it Proof.} Since $V(P)$ is a tree, then by Lemma~\ref{l0}, $V_1(P)=\emptyset$, $V_2(P)\not=\emptyset$ and, $V(P)=V_2(P)\cup V_3(P).$

To prove the lemma, we try to construct a probability matrix $Q$ such that
\begin{itemize}
  \item $Q=P(\sigma,\pi)\in {\cal C}(\sigma^{-1} {\bf p},\pi^{-1} {\bf q})$, for some permutation pair $(\sigma,\pi)\in \Sigma_m\times\Sigma_n$;
  \item $|V(Q)|=m+n-1$.
\end{itemize}

To define the matrix $Q=(q_{i,j})_{m\times n}$, firstly, for a fixed $(i_0,j_0)\in V_2(P)$, without loss of generality, suppose that $p_{i_0,j_0}$ be the unique positive entry in the $i_0$-th row of $P$. We define a {\it directed} path $\gamma$ in $V(Q)$ as follows.
\begin{itemize}
  \item Let $q_{1,1}=p_{i_0,j_0}$, denote $(i(0),j(0))=(1,1)$.
  \item Write $\{l\not=i_0:(l,j_0)\in V_2(P)\}=\{l_k:1\leq k \leq s_1\}$, $s_1\geq 0$, such that $p_{l_k,j_0}$ decreases in $k$, let $q_{i(0)+k,j(0)}=p_{l_k,j_0}$, $1\leq k\leq s_1$. Define $U_1:=\{(l,j_0):(l,j_0)\in V_3(P)\}$, let $(i_1,j_0)$ be the element in $U_1$ such that $p_{i_1,j_0}=\max\{p_{i,j}:(i,j)\in U_1\}$ (note that in this way, we define $i_1$), denote $i(1):=i(0)+s_1+1$, and let $q_{i(1),j(0)}=p_{i_1,j_0}$.
  \item Write $\{t\not=j_0:(i_1,t)\in V_2(P)\}=\{t_k:1\leq k \leq s_2\}$ such that $p_{i_1,t_k}$ decreases in $k$, let $q_{i(1),j(0)+k}=p_{i_1,t_k}$, $1\leq k\leq s_2$. Define $U_2:=\{(i_1,t):(i_1,t)\in V_3(P)\setminus\{(i_1,j_0)\}\}$, let $(i_1,j_1)$ be the element in $U_2$ such that $p_{i_1,j_1}=\max\{p_{i,j}:(i,j)\in U_2\}$, denote $j(1):=j(0)+s_2+1$, and let $q_{i(1),j(1)}=p_{i_1,j_1}$.

  \item For $\xi\geq 3$. In the case when $\xi=2\zeta-1$, write $\{l\not=i_{\zeta-1}:(l,j_{\zeta-1})\in V_2(P)\}=\{l_k:1\leq k \leq s_{2\zeta-1}\}$ such that $p_{l_k,j_{\zeta-1}}$ decreases in $k$, let $q_{i(\zeta-1)+k,j(\zeta-1)}=p_{l_k,i_{\zeta-1}}$, $1\leq k\leq s_{2\zeta-1}$. Define $U_\xi:=\{(l,j_{\zeta-1}):(l,j_{\zeta-1})\in V_3(P)\}$, let $(i_\zeta,j_{\zeta-1})$ be the element in $U_\xi$ such that $p_{i_\zeta,j_{\zeta-1}}=\max\{p_{i,j}:(i,j)\in U_\xi\}$, denote $i(\zeta):=i(\zeta-1)+s_{2\zeta-1}+1$, and let $q_{i(\zeta),j(\zeta-1)}=p_{i_\zeta,j_{\zeta-1}}$.

      In the case when $\xi=2\zeta$, write $\{t\not=j_{\zeta-1}:(i_\zeta,t)\in V_2(P)\}=\{t_k:1\leq k \leq s_{2\zeta}\}$ such that $p_{i_\zeta,t_k}$ decreases in $k$, let $q_{i(\zeta),j(\zeta-1)+k}=p_{i_\zeta,t_k}$, $1\leq k\leq s_{2\zeta}$. Define $U_\xi:=\{(i_\zeta,t):(i_\zeta,t)\in V_3(P)\}$, let $(i_\zeta,j_{\zeta})$ be the element in $U_\xi$ such that $p_{i_\zeta,j_{\zeta}}=\max\{p_{i,j}:(i,j)\in U_\xi\}$, denote $j(\zeta):=j(\zeta-1)+s_{2\zeta}+1$, and let $q_{i(\zeta),j(\zeta)}=p_{i_\zeta,j_{\zeta}}$.

      \item Repeat the above procedure for $\xi\geq 1$ until $U_{\xi}=\emptyset$. Let $\xi_0=\min\{\xi\geq 1:|U_\xi|=0\}$. When $\xi_0=2\zeta_0-1$, $\gamma$ is the directed path in $Q$ from $(i(0),j(0))$ to $(i(\zeta_0-1)+s_{\xi_0},j(\zeta_0-1))$; when $\xi_0=2\zeta_0$, $\gamma$ is the directed path in $Q$ from $(i(0),j(0))$ to $(i(\zeta_0),j(\zeta_0-1)+s_{\xi_0})$.

\end{itemize}

If $|U_k|=1$ for all $1\leq k\leq \xi_0-1$, all points in $V_3(P)$ and then all points in $V_2(P)$ are used in the definition of $\gamma$, so, $|V(P)|=|\gamma|$. On the other hand, since $V(P)$ is complete, then $\gamma$ forms a {\it continuous directed} path from $(1,1)$ to $(m,n)$, this implies $|\gamma|=m+n-1$. For any $(i,j)\notin \gamma$, let $q_{i,j}=0$, thus, we obtain $Q$ as required.

Otherwise, write $\gamma_0=\gamma$ and let $V(\gamma_0)=\{(i,j)\in V(P):p_{i,j}$ is used in the definition of $\gamma_0\}$.  We re-define $U_k:=U_k\setminus V(\gamma_0)$, for any $1\leq k\leq \xi_0$, let $k_0:=\max\{k<\xi_0:|U_k|>1\}$. Without loss of generality, suppose that $k_0=2z_0$ is even. By the definition of $U_{k_0}$, $U_{k_0}\subset V_3(P)$ and for any $(i,j)\in U_{k_0}$, $i=i_{z_0}$ is defined in the definition of $\gamma$. Let $(i_{z_0},j_{\zeta_0})\in U_{k_0}$ such that $p_{i_{z_0},j_{\zeta_0}}=\max\{p_{i,j}:(i,j)\in U_{k_0}\}$.
Without of loss of generality, we assume $\xi_0=2\zeta_0-1$, recall that in this case the end vertex of $\gamma_0$ is $\(i(\zeta_0-1)+s_{2\zeta_0-1},j(\zeta_0-1)\)$. Denote $(i(\zeta_0),j(\zeta_0)):=\(i(\zeta_0-1)+s_{2\zeta_0-1},j(\zeta_0-1)+1\)$, $s_{2\zeta_0}:=0$ and $U_{\xi_0+1}=U_{2\zeta_0}=\emptyset$.

Now, similar to $\gamma_0$, we define another directed path $\gamma_1$ in $V(Q)$ from $(i(z_0),j(\zeta_0))$ as follows.
\begin{itemize}
  \item Let $q_{i(z_0),j(\zeta_0)}=p_{i_{z_0},j_{\zeta_0}}$.
  \item Write $\{l\not=i_{z_0}:(l,j_{\zeta_0})\in V_2(P)\}=\{l_k:1\leq k \leq s_{2\zeta_0+1}\}$, such that $p_{l_k,j_{\zeta_0}}$ decreases in $k$, let $q_{i(\zeta_0)+k,j(\zeta_0)}=p_{l_k,j_{\zeta_0}}$, $1\leq k\leq s_{2\zeta_0+1}$. Define $U_{\xi_0+2}=U_{2\zeta_0+1}:=\{(l,j_{\zeta_0}):(l,j_{\zeta_0})\in V_3(P)\setminus{(i_{z_0},j_{\zeta_0})}\}$, let $(i_{\zeta_0+1},j_{\zeta_0})\in U_{\xi_0+2}$ such that $p_{i_{\zeta_0+1},j_{\zeta_0}}=\max\{p_{i,j}:(i,j)\in U_{\xi_0+2}\}$, denote $i(\zeta_0+1)=i(\zeta_0)+s_{2\zeta_0+1}+1$, let $q_{i(\zeta_0+1),j(\zeta_0)}=p_{i_{\zeta_0+1},j_{\zeta_0}}$.
\item Write $\{t\not=j_{\zeta_0}:(i_{\zeta_0+1},t)\in V_2(P)\}=\{t_k:1\leq k \leq s_{2(\zeta_0+1)}\}$, such that $p_{i_{\zeta_0+1},t_k}$ decreases in $k$, let $q_{i(\zeta_0+1),j(\zeta_0)+k}=p_{i_{\zeta_0+1},t_k}$, $1\leq k\leq s_{2(\zeta_0+1)}$. Define $U_{\xi_0+3}=U_{2(\zeta_0+1)}:=\{(i_{\zeta_0+1},t):(i_{\zeta_0+1},t)\in V_3(P)\setminus{(i_{\zeta_0+1},j_{\zeta_0})}\}$, let $(i_{\zeta_0+1},j_{\zeta_0+1})\in U_{\xi_0+3}$ such that $p_{i_{\zeta_0+1},j_{\zeta_0+1}}=\max\{p_{i,j}:(i,j)\in U_{\xi_0+3}\}$, denote $j(\zeta_0+1)=j(\zeta_0)+s_{2(\zeta_0+1)}+1$, let $q_{i(\zeta_0+1),j(\zeta_0+1)}=p_{i_{\zeta_0+1},j_{\zeta_0+1}}$.
\item $\cdots\cdots$
\item Repeat the above procedure for $\xi\geq \xi_0+2$ until $U_{\xi}=\emptyset$, and let $\xi_1=\min\{k\geq\xi_0+2:U_{\xi}=\emptyset\}$. When $\xi_1=2\zeta_1-1$, $\gamma_1$ is the directed path in $Q$ from $(i(z_0),j(\zeta_0))$ to $(i(\zeta_0)+1,j(\zeta_0))$ and then to $(i(\zeta_1-1)+s_{\xi_1},j(\zeta_1-1))$; when $\xi_1=2\zeta_1$, $\gamma_1$ is the directed path in $Q$ from $(i(z_0),j(\zeta_0))$ to $(i(\zeta_0)+1,j(\zeta_0))$ and then to $(i(\zeta_1),j(\zeta_1-1)+s_{\xi_1})$.
\end{itemize}

If $\cup_{k=1}^{\xi_1-1}U_k\setminus{V\(\gamma_0\cup\gamma_1\)}=\emptyset$, where $V(\gamma_0\cup\gamma_1)$, together with the following $V(\cup^\tau_{k=0}\gamma_\tau)$, is same defined as $V(\gamma_0)$, then $|V(P)|=|\gamma_0\cup\gamma_1|$ and $\(\gamma_0\cup\gamma_1\cup\{(i(\zeta_0),j(\zeta_0))\}\)\setminus\{(i(z_0),j(\zeta_0))\}$ forms a continuous directed path from $(1,1)$ to $(m,n)$. Thus $|V(P)|=|\gamma_0\cup\gamma_1|=m+n-1$. For any $(i,j)\notin \gamma_0\cup \gamma_1$, let $q_{i,j}=0$, we obtain $Q$ as required.

If $\cup_{k=1}^{\xi_1-1}U_k\setminus{V\(\gamma_0\cup\gamma_1\)}\not=\emptyset$, re-define $U_k=U_k\setminus V(\gamma_0\cup\gamma_1)$ for any $1\leq k\leq \xi_1$. Let $k_1:=\max\{k<\xi_1:|U_k|>1\}$. Without loss of generality, suppose that $\xi_1=2\zeta_1$ and $k_1=2z_1-1$. By the definition of $U_{k_1}$, $U_{k_1}\subset V_3(P)$ and for any $(i,j)\in U_{k_1}$, $j=j_{z_1-1}$ is defined in the definition of $\gamma_0$ and $\gamma_1$. Let $(i_{\zeta_1+1},j_{z_1-1})\in U_{k_1}$ such that $p_{i_{\zeta_1+1},j_{z_1-1}}=\max\{p_{i,j}:(i,j)\in U_{k_1}\}$. Recall that in this case the end vertex of $\gamma_1$ is $\(i(\zeta_1),j(\zeta_1-1)+s_{2\zeta_1}\)$.
Denote $(i(\zeta_1+1),j(\zeta_1)):=\(i(\zeta_1)+1,j(\zeta_1-1)+s_{2\zeta_1}\)$, $s_{2\zeta_1+1}:=0$ and $U_{\xi+1}=U_{2\zeta_1+1}:=\emptyset$. Similar to $\gamma_0,\gamma_1$, we define a directed path $\gamma_2$ in $V(Q)$ from $(i(\zeta_1+1),j(z_1-1))$ by defining $q_{i(\zeta_1+1),j(z_1-1)}=p_{i_{\zeta_1+1},j_{z_1-1}}$, $\cdots\cdots$.

$\cdots\cdots$

We stop until we obtain a directed path $\gamma_\tau$, which ends at the vertex $(m,n)$. For any $1\leq k\leq \tau$, denote by $u_k$ the beginning point of $\gamma_k$, $w_k$ the second point of $\gamma_k$ and $v_k$ the vertex in the interval between $u_k$ and $w_k$ such that the Euclidean distance between $v_k$ and $w_k$ is 1. For example, in our construction, $u_1=(i(z_0),j(\zeta_0))$, $u_2=(i(\zeta_1+1),j(z_1-1))$; $v_1=(i(\zeta_0),j(\zeta_0))$, $v_2=(i(\zeta_1+1),j(\zeta_1))$. Let $\bar\gamma_k=\(\gamma_k\cup\{v_k\}\)\setminus\{u_k\}$ for any $1\leq k\leq \tau$, then $\bar \gamma_k$ forms a continuous directed path beginning at $v_k$. Denote $\bar\gamma_0=\gamma_0$, thus $\cup^\tau_{k=0}\bar\gamma_\tau$ forms a continuous directed path from $(1,1)$ to $(m,n)$ and $|\cup^\tau_{k=0}\gamma_\tau|=|\cup^\tau_{k=0}\bar\gamma_\tau|$.

\begin{figure}
\unitlength=0.35mm
\begin{picture}(300,300)(-90,0)
\multiput(20,0)(10,0){26}{\circle{1}}\multiput(260,0)(10,0){2}{\circle*{3}}\put(247.5,-2.5){$\star$}
\multiput(20,10)(10,0){26}{\circle{1}}\put(250,10){\circle*{3}}
\multiput(20,20)(10,0){26}{\circle{1}}\put(17.5,17.5){$\star$}\put(240,20){\circle*{3}}\put(247.5,17.5){$\star$}
\multiput(20,30)(10,0){26}{\circle{1}}\multiput(220,30)(10,0){2}{\circle*{3}}\put(207.5,27.5){$\star$}
\multiput(20,40)(10,0){26}{\circle{1}}\put(210,40){\circle*{3}}
\multiput(20,50)(10,0){26}{\circle{1}}\put(207.5,47.5){$\star$}\put(187.5,47.5){$\star$}\put(200,50){\circle*{3}}
\multiput(20,60)(10,0){26}{\circle{1}}\put(190,60){\circle*{3}}
\multiput(20,70)(10,0){26}{\circle{1}}\put(180,70){\circle*{3}}\put(167.5,67.5){$\star$}
\multiput(20,80)(10,0){26}{\circle{1}}\put(170,80){\circle*{3}}
\multiput(20,90)(10,0){26}{\circle{1}}\put(157.5,87.5){$\star$}\put(167.5,87.5){$\star$}
\multiput(20,100)(10,0){26}{\circle{1}}\put(160,100){\circle*{3}}
\multiput(20,110)(10,0){26}{\circle{1}}\put(150,110){\circle*{3}}\put(157.5,107.5){$\star$}\put(37.5,107.5){$\star$}
\multiput(20,120)(10,0){26}{\circle{1}}\multiput(130,120)(10,0){2}{\circle*{3}}\put(117.5,117.5){$\star$}
\multiput(20,130)(10,0){26}{\circle{1}}\put(120,130){\circle*{3}}
\multiput(20,140)(10,0){26}{\circle{1}}\multiput(100,140)(10,0){2}{\circle*{3}}\multiput(87.5,137.5)(30,0){2}{$\star$}
\multiput(20,150)(10,0){26}{\circle{1}}\put(90,150){\circle*{3}}
\multiput(20,160)(10,0){26}{\circle{1}}\put(80,160){\circle*{3}}
\multiput(20,170)(10,0){26}{\circle{1}}\put(80,170){\circle*{3}}
\multiput(20,180)(10,0){26}{\circle{1}}\put(77.5,177.5){$\star$}\put(57.5,177.5){$\star$}\put(70,180){\circle*{3}}
\multiput(20,190)(10,0){26}{\circle{1}}\put(60,190){\circle*{3}}
\multiput(20,200)(10,0){26}{\circle{1}}\put(60,200){\circle*{3}}
\multiput(20,210)(10,0){26}{\circle{1}}\put(57.5,207.5){$\star$}\put(37.5,207.5){$\star$}\put(50,210){\circle*{3}}\put(87.5,207.5){$\star$}
\multiput(20,220)(10,0){26}{\circle{1}}\put(40,220){\circle*{3}}
\multiput(20,230)(10,0){26}{\circle{1}}\put(40,230){\circle*{3}}
\multiput(20,240)(10,0){26}{\circle{1}}\put(37.5,237.5){$\star$}\put(187.5,237.5){$\star$}\put(30,240){\circle*{3}}\put(17.5,237.5){$\star$}
\multiput(20,250)(10,0){26}{\circle{1}}\put(20,250){\circle*{3}}
\multiput(20,260)(10,0){26}{\circle{1}}\put(20,260){\circle*{3}}
\put(18,262){\line(1,0){64}}\put(18,158){\line(1,0){64}}\put(18,158){\line(0,1){104}}\put(82,158){\line(0,1){104}}
\put(88,118){\line(1,0){54}}\put(88,212){\line(1,0){54}}\put(88,118){\line(0,1){94}}\put(142,118){\line(0,1){94}}
\put(38,112){\line(1,0){144}}\put(38,68){\line(1,0){144}}\put(38,68){\line(0,1){44}}\put(182,68){\line(0,1){44}}
\put(18,-2){\line(1,0){254}}\put(18,22){\line(1,0){254}}\put(18,-2){\line(0,1){24}}\put(272,-2){\line(0,1){24}}
\put(188,28){\line(0,1){214}}\put(232,28){\line(0,1){214}}\put(188,28){\line(1,0){44}}\put(188,242){\line(1,0){44}}
\put(70,254){$\gamma_0$}\put(130,204){$\gamma_1$}\put(40,74){$\gamma_2$}\put(220,234){$\gamma_3$}\put(20,2){$\gamma_4$}
\put(90,160){\color{red}\circle*{3}}\put(140,110){\color{red}\circle*{3}}\put(190,70){\color{red}\circle*{3}}\put(230,20){\color{red}\circle*{3}}
\put(90,162){$v_1$} \put(90,202){$u_1$}\put(133,103){$v_2$}\put(40,103){$u_2$}\put(190,62){$v_3$}\put(190,232){$u_3$}\put(223,13){$v_4$}\put(20,13){$u_4$}
\put(186,28){\color{red}\line(1,0){48}}\put(186,28){\color{red}\line(0,1){44}}\put(186,72){\color{red}\line(1,0){48}}\put(234,28){\color{red}\line(0,1){44}}
\put(220,62){$\bar\gamma_3$}
\put(-6,258){\footnotesize(1,1)}\put(274,-3){(m,n)=\footnotesize(27,26)}
\end{picture}
\vskip 0mm
\begin{center}
\begin{minipage}{13.5cm}

\caption{\it\small An example to show the structure of $Q$:
a) $\bullet$, $\star$ is the new position of a point in $V_2(P)$, $V_3(P)$ respectively, $\circ$'s are zeros. $|V(Q)|=|V(P)|$ and $V(Q)=\cup_{k=0}^{\tau}\gamma_k$, $\tau=4$. For each $1\leq k\leq 4$, $u_k$ is the beginning vertex of $\gamma_k$, $v_k$ ($\notin V(Q)$) is the beginning vertex of $\bar\gamma_k$, and $\gamma_0\cup\bar\gamma_1\cup\cdots\cup\bar\gamma_4$ forms a continuous directed path from $(1,1)$ to $(m,n)=(27,26)$. b) In the definition of $\gamma_0$, one has $\xi_0=2\zeta_0-1=7$, $k_0=2z_0=4$, $s_1=s_2=s_4=s_6=1$, $s_3=s_5=s_7=2$, $i(1)=3,\ i(2)=6,\ i(3)=9$, $j(1)=3,\ j(2)=5,\ j(3)=7$. Before we define $\gamma_1$, we define $(i(4),j(4))=(11,8)=v_1$. c) The minimum-entropy coupling possesses nice local features, for example, by Theorem~\ref{t2}, the subpath $\gamma_0'$ $=\{(1,1),(3,1),(3,3),(6,3),(6,5),(9,5),(9,7),(11,7)\}$, which forms the skeleton of $\gamma_0$, behaves supper-Fibonacci, i.e. $q_{1,1}+q_{3,1}\leq q_{3,3}$, $q_{3,1}+q_{3,3}\leq q_{6,3}$,\ldots, $q_{9,5}+q_{9,7}\leq q_{11,7}$.} \label{f1}
\end{minipage}
\end{center}
\end{figure}
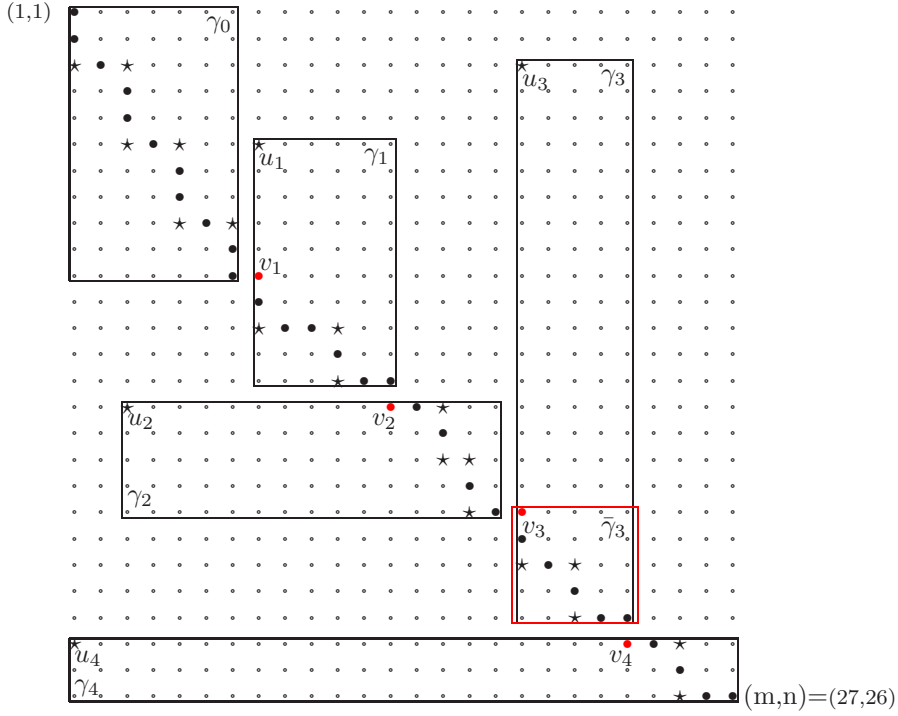

Now, $V(\cup^\tau_{k=0}\gamma_\tau)=V(P)$, and $|V(P)|=|\cup^\tau_{k=0}\gamma_\tau|$=$|\cup^\tau_{k=0}\bar\gamma_\tau|$=$m+n-1$. Finally, for any $(i,j)\notin \cup^\tau_{k=0}\gamma_\tau$, define $q_{i,j}=0$, then we obtain $Q$ as required. For an illustration of the structure of $Q$, see Figure~\ref{f1}.\QED

{\it Proof of Theorem~\ref{t3}.} By Theorem~\ref{t2}, $V(P)$ is a complete forest, suppose that the number of connected components of the forest is $k$, $k\geq 1$. Then there exists some permutation pair $(\sigma,\pi)\in\Sigma_m\times\Sigma_n$ such that $P(\sigma,\pi)\in {\cal C}(\sigma^{-1}{\bf p},\pi^{-1}{\bf q})$ has the following block structure

\be\label{2.6}P(\sigma,\pi)=\(\ba{llll}P_{1}&0&\ldots& 0\\
0&P_{2}&\ldots&0\\
\vdots&\vdots&\ddots&\vdots\\
0&0&\ldots&P_{k}
\ea\)\ee
Where $P_{l}$ is a $m_l\times n_l$ submatrix with $1\leq m_l\leq m, 1\leq n_l\leq n$, $1\leq l\leq k$ and $\sum_{l=1}^km_l=m$, $\sum_{l=1}^kn_l=n$. Furthermore, $V(P_{l})$ is a complete tree in $\{M_{l-1}+1,\ldots, M_{l-1}+m_l\}\times\{N_{l-1}+1,\ldots, N_{l-1}+n_l\}$, where $M_l=\sum_{i=1}^lm_i$, $N_l=\sum_{i=1}^ln_i$, $1\leq l\leq k$. The $0$'s are the corresponding zero matrixes.

By Lemma~\ref{l3}, $|V(P_{l})|=n_l+m_l-1$ and then $$|V(P)|=|V(P(\sigma,\pi))|=\sum_{l=1}^k|V(P_{l})|=m+n-k.$$
Finally, by the block structure of $P(\sigma,\pi)$ and the definition of $\kappa(\bf p,q)$, it holds $1\leq k\leq \kappa(\bf p,q)$, then the theorem follows.\QED

{\it Proof of Theorem~\ref{t1}}. First of all, for any $P=(p_{i,j})_{m\times n}\in{\cal C}_e({\bf p},{\bf q})$, we claim that $V(P)$ is a complete forest. The following proof is based on a private discussion with Professor Yu Lei. In fact, if there is a circuit $\gamma=\{v_0,v_1,\ldots, v_s=v_0\}$ in $V(P)$ ($s$ is even and $\geq 4$), take $0<\epsilon<\min\{p_{v_i}:1\leq i\leq s\}$ and define $P'=(p'_{i,j})_{m\times n}$, $P''=(p''_{i,j})_{m\times n}$ as the following:
$$ p'_{i,j}=\left\{\ba{ll}p_{i,j},&{\rm if}\ (i,j)\notin \gamma;\\
p_{i,j}+\epsilon,& {\rm if}\ (i,j)=v_k,\ \ k\ {\rm is\ even;}\\
p_{i,j}-\epsilon,& {\rm if}\ (i,j)=v_k,\ \ k\ {\rm is\ odd,}\ea\right.{\rm and}\ p''_{i,j}=\left\{\ba{ll}p_{i,j},&{\rm if}\ (i,j)\notin \gamma;\\
p_{i,j}+\epsilon,& {\rm if}\ (i,j)=v_k,\ \ k\ {\rm is\ odd;}\\
p_{i,j}-\epsilon,& {\rm if}\ (i,j)=v_k,\ \ k\ {\rm is\ even.}\ea\right.
$$
Then $P',P''\in{\cal C}({\bf p},{\bf q})$ and $P=\frac 12 P'+\frac 12 P''$, a contradiction.

By the proof of Theorem \ref{t3}, we have
$$ P\in \Dp\bigcup_{l=m+n-\kappa({\bf p},{\bf q})}^{m+n-1}{\hC}_l({\bf p},{\bf q}),$$
where ${\hC }_l=\{P\in{\cal C}({\bf p},{\bf q}): {\rm for\ some\ complete\ forest\ }F\ {\rm with\ }|F|=l,\ V(P)=F\}$.

Since any complete forest $F$ is a subgraph of some tree $T\in {\cal T}$, we have
$$\bigcup_{l=m+n-\kappa({\bf p},{\bf q})}^{m+n-1}{\hC }_l({\bf p},{\bf q})=\hC ({\bf p},{\bf q}).$$ Thus ${\cal C}_e({\bf p},{\bf q})\subset {\hC}({\bf p},{\bf q})$.

Second, for any $P\in {\hC}({\bf p},{\bf q})$, if $P$ is not an extreme point, then there exists $P_1, P_2, \ldots, P_l\in {\cal C}({\bf p},{\bf q})$ and $\lambda_1,\lambda_2,\ldots,\lambda_l\in(0,1)$, $l\geq 2$, such that $$P=\sum_{i=1}^l\lambda_iP_i.$$
Then $V(P)=\cup_{i=1}^lV(P_i)$ and we have
$$V(P_i)\subset V(P)\subset T,\ \forall\ 1\leq i\leq l,$$ for some $T\in {\cal T}$. By Proposition~\ref{p1}, there exists at most one $P\in{\cal C}({\bf p},{\bf q})$ such that $P$ is consistent with $T$, one has $P_1=P_2=\ldots=P_l=P$. So, $P\in {\cal C}_e({\bf p},{\bf q})$ and
$ {\hC}({\bf p},{\bf q})\subset {\cal C}_e({\bf p},{\bf q})$.
\QED
Actually, by the above arguments, a coupling $P\in{\cal C}({\bf p},{\bf p })$ which can not be further optimized by Lemma 1 is an extreme point; if we optimize such a $P$ to $P'$ by Lemma 2 or Lemma~\ref{l1}, then $P'$ is another extreme point such that $H(P')<H(P)$. Note that the so-called {\it greedy coupling} $P$ provided by the greedy algorithm, which is first posed in \cite{KDVH} and then developed in \cite{KDVH2} and \cite{CK,R} etc, possesses the forest structure and is an extreme point of ${\cal C}({\bf p},{\bf p })$.

Finally, we have the following corollary.

\bco\label{cor2} For any function $\Psi$ on ${\cal C}({\bf p},{\bf q})$, if for any minimal value point $\tilde P$ of $\Psi$, $V(\tilde P)$ is a forest, then $\Psi$ takes its minimal value in ${\cal C}_e({\bf p},{\bf q})={\hC}({\bf p},{\bf q})$. \eco

\vskip 5mm
\section{The algorithm via an algebraic argument}
\renewcommand{\theequation}{3.\arabic{equation}}
\setcounter{equation}{0}

Let $\hS$ be the collection of subsets of $[mn]$ with cardinality $m+n-1$. For any $S\in\hS$, enumerate $S=\{i_1,i_2,\ldots,i_{m+n-1}\}$ such that $i_1<i_2<\ldots<i_{m+n-1}$. For any $S=\{i_1,i_2,\ldots,i_{m+n-1}\}$, $S'=\{j_1,j_2,\ldots,j_{m+n-1}\}\in \hS$, as usual, we call $S\prec S'$ in lexicographic order if and only if for some $0\leq k_0\leq m+n-2$, $i_k=j_k$ for $1\leq k\leq k_0$ and $i_{k_0+1}<j_{k_0+1}$. Let $S_k$ be the $k$-th element of $\hS $ in lexicographic order for $1\leq k\leq C(m,n)=\binom{mn}{m+n-1}$.

For any $1\leq i\leq mn$, define $\phi(i):=(t,r)\in[m]\times[n]$, where $(t,r)$ is the unique element in $[m]\times[n]$ such that
\be\label{3.1}i=(t-1)n+r.\ee

Denote by ${\cal V}$ the collection of subsets of $[m]\times[n]$ with cardinality $m+n-1$. Let $V_k=\phi(S_k):=\{\phi(i):i\in S_k\}$, $1\leq k\leq C(m,n)$. Clearly, $\phi$ is a bijection between $\hS$ and $\cal V$, here with a little abuse of notation, we call $V_k$ the $k$-th element in $\cal V$ in lexicographic order.

\bd\label{d10} For any $S=\{i_1,i_2,\ldots,i_{m+n-1}\}\in \hS$, let $V=\phi(S)\in{\cal V}$. We call the $(m+n-1)\times(m+n-1)$ matrix $\fA=\fA(S)=\fA(V)=(a_{s,k})$ defined below the \textbf{structure matrix} of $S$ or $V$. For any $1\leq k\leq m+n-1$, if $\phi(i_k)=(t,r)$, then
\be\label{3.2} a_{s,k}:=\left\{\ba{ll}1,& {\rm if }\ s=t<m;\\
1,&{\rm if}\ s=m;\\
1,&{\rm if}\ s=m+r<m+n;\\
0,&{\rm else.}\ea\right.\ee
\ed
For any $c>0$, let ${\cal P}^+_{m,c}=\{{\bf p}=(p_1,\ldots,p_m):\sum_{k=1}^mp_k=c,\ p_k>0,\ 1\leq k\leq m\}$. For any ${\bf p}\in{\cal P}^+_{m,c},\ {\bf q}\in{\cal P}^+_{n,c}$, let ${\cal M}_c({\bf p,q})$ be the collection of $m\times n$ matrix $B=(b_{i,j})$ such that
$$\sum_{i=1}^mb_{i,j}=q_j,\ \sum_{j=1}^nb_{i,j}=p_i,\ {\rm for \ all\ }1\leq i\leq m,\ i\leq j\leq n.$$
For any $B\in {\cal M}_c({\bf p,q})$, let $V(B)=\{(i,j):b_{i,j}\not=0\}$. In the case of  $c=1$, writing ${\cal M}({\bf p,q})={\cal M}_1({\bf p,q})$, one has ${\cal C}({\bf p,q})\subset {\cal M}({\bf p},{\bf q})$.

For any ${\bf p}\in{\cal P}^+_{m,c},\ {\bf q}\in{\cal P}^+_{n,c}$, let $y_c({\bf p,q})=(p_1,\ldots,p_{m-1},c,q_1,\ldots,q_{n-1})^T$, where $(\cdot)^T$ means the transpose of $(\cdot)$ and $y_c({\bf p,q})$ is a column vector in $\R^{m+n-1}$.

\bl\label{l4} For any $V\in{\cal V}$, suppose $\phi^{-1}(V)=S=\{i_1,i_2,\ldots,i_{m+n-1}\}\in \hS$. Then there exists $B\in {\cal M}_c({\bf p,q})$ such that $V(B)\subset V$ if and only if the following system of linear equations has a solution $x=(x_1,x_2,\ldots,x_{m+n-1})^T\in\R^{m+n-1}$:
\be\label{3.3} \fA(V)x=y_c(\bf p,q).\ee In particular, the solution $x$ and the matrix $B$ are determined from each other in the following fashion:
$b_{\phi(i_k)}=x_k$, $1\leq k\leq m+n-1$;
$b_{t,r}=0$ else.
\el
{\it Proof.} The lemma follows straightforwardly from the definition of the structure matrix.\QED

By introducing the concept of structure matrix, we obtain the following criteria theorem for trees in ${\cal V}$.
\bt For any $V\in{\cal V}$, $V$ is a tree if and only if $\det(\fA(V))\not=0$, i.e. $\fA(V)$ is reversible. Where $\fA(V)$ is the structure matrix of $V$ defined in (\ref{3.2}).\et

{\it Proof.}
Let's begin with the necessary part of the proof. If $V\in{\cal V}$ is a tree, first of all, by the proof of Proposition~\ref{p1}, in any case there exists a unique $B\in {\cal M}(\bf p,q)$ such that $V(B)\subset V$. Furthermore, for any $c>0$, for any ${\bf p}\in{\cal P}^+_{m,c},\ {\bf q}\in{\cal P}^+_{n,c}$, there exists a unique $B\in {\cal M}_c(\bf p,q)$ such that $V(B)\subset V$. By Lemma~\ref{l4}, the latter is equivalent to the fact that there exists a unique $x=(x_1,x_2,\ldots,x_{m+n-1})^T\in \R^{m+n-1}$ such that
$\fA(V)x=y_c(\bf p,q).$ Clearly, the unique solution $x\not=(0,0,\ldots,0)^T$.

Now, fix ${\bf p}\in{\cal P}^+_m$, ${\bf q}\in{\cal P}^+_n$ arbitrarily, let $y({\bf p,q})=(p_1,\ldots,p_{m-1},1,q_1,\ldots,q_{n-1})^T$. Take $\epsilon=\epsilon({\bf p,q})>0$ small enough such that, for any $y=(y_1,y_2,\ldots,$ $y_{m+n-1})^T\in {\cal B}(y({\bf p,q}),\epsilon)$, $y_k>0$ for all $k$ and $y_m-\sum_{k=1}^{m-1}y_k>0$, $y_m-\sum_{k=1}^{n-1}y_{m+k}>0$. Where ${\cal B}(y({\bf p,q}),\epsilon)\subset \R^{m+n-1}$ is the ball with radius $\epsilon$ centered at $y(\bf p,q)$.

For any $y\in {\cal B}(y(\bf p,q),\epsilon)$, take $c=y_m$, ${\bf p}=(y_1,y_2,\ldots,y_{m-1}, y_m-\sum_{k=1}^{m-1}y_k)$ and ${\bf q}=(y_{m+1},y_{m+2},\ldots,$ $y_{m+n-1}, y_m-\sum_{k=1}^{n-1}y_{m+k})$, then ${\bf p}\in{\cal P}^+_{m,c},\ {\bf q}\in{\cal P}^+_{n,c}$ and $y=y_c(\bf p,q)$. By the arguments in the first paragraph of the proof, for $y_c({\bf p,q}) (=y)$, equation (\ref{3.3}) has a unique solution $x\in\R^{m+n-1}$ and $x\not=(0,0,\ldots,0)^T$. This implies that the $(m+n-1)$-dimensional ball ${\cal B}(y(\bf p,q),\epsilon)$ is contained in the linear space spanned by the column vectors of $\fA(V)$, hence the column vectors of $\fA(V)$ are linearly independent and $\det(\fA(V))\not=0$.

For the sufficiency part of the proof, we assume $V$ is not a tree. To show $\det(\fA(V))=0$, it suffices to prove that there exists ${\bf p}\in {\cal P}^+_m$, ${\bf q}\in {\cal P}^+_n$, such that equation (\ref{3.3}) has no solution. By Lemma~\ref{l4}, this is equivalent that there is no $B\in {\cal M}(\bf p,q)$ such that $V(B)\subset V$.

To this end, take ${\bf p}\in {\cal P}^+_m$, ${\bf q}\in {\cal P}^+_n$ such that $\kappa(\bf p,q)=1$, then, same as Proposition~\ref{p-}, for any $B\in {\cal M}(\bf p,q)$, $V(B)$ is complete and connected. By Proposition~\ref{p0}, iii), this implies $|V(B)|\geq m+n-1(=|V|)$. Now, if there exists some $B\in {\cal M}(\bf p,q)$ such that $V(B)\subset V$, then $V(B)=V$ and $V(B)$ is not a tree, a contradiction to Proposition~\ref{p0}, ii) appears. Thus $\det(\fA(V))=0$.\QED

\bt\label{t4'} For any ${\bf p}\in {\cal P}^+_m,$ ${\bf q}\in{\cal P}^+_n$, let $y({\bf p,q})=(p_1,\ldots,p_{m-1},1,q_1,\ldots,q_{n-1})^T$. For any $1\leq k\leq C(m,n)$, denote $\fA_k:=\fA(S_k)=\fA(V_k)$, the structure matrix of $S_k$ or $V_k$. If $\det(\fA_k)\not=0$ and $\fA_k^{-1}y({\bf p,q})\in {\cal P}_{m+n-1}$, denote by $P_k$ the coupling determined by $\fA_k^{-1}y({\bf p,q})$ as in the statement of Lemma~\ref{l4}. Then
\be\label{3.40}{\cal C}_e({\bf p},{\bf q})={\hC}({\bf p},{\bf q})=\left\{P_k:\det(\fA_k)\not=0\ {\rm  and}\ \fA_k^{-1}y({\bf p,q})\in {\cal P}_{m+n-1},\ 1\leq k\leq C(m,n)\right\}.\ee

For $H$, a strictly concave function as Shannon entropy, R\'enyi entropy or Tsallas entropy on ${\cal C}({\bf p},{\bf q})$, define
\be\label{3.4}H_k:=\left\{\ba{ll}H(P_k),& \ {\rm if} \det(\fA_k)\not=0\ {\rm  and}\ \fA_k^{-1}y({\bf p,q})\in {\cal P}_{m+n-1};\\[3mm]\infty, &\ {\rm otherwise.}\ea\right.\ee
Then \be \inf_{P\in{\cal C}({\bf p},{\bf q})}H(P)=\min_{P\in \hC({\bf p,q})}H(P)=\min\left\{H_k:1\leq k\leq C(m,n)\right\}.\ee \et

The following is the algorithm to calculate the minimal joint entropy and the corresponding minimum-entropy couplings.

\noindent\hrulefill

\noindent{\textbf{Algorithm}}: The Min Entropy Coupling Algorithm
\vskip-2mm
\noindent\hrulefill

\noindent MIN-ENTROPY-COUPLING~$({\bf p},{\bf q})$

\noindent{\textbf{Input:}} probability distributions~${\bf p}=(p_1,p_2,\ldots, p_m)$ and ${\bf q}=(q_1,q_2,\ldots, q_n)$, $\mathscr{S}=\left\{S_k: 1\leq k\leq C(m,n)\right\}$ be the collection of subsets of~$[mn]$~with cardinality~$m+n-1$.

\noindent{\textbf{Output:}} An~$n\times n$ matrix~$P =(p_{i,j})$ s.t. $\sum_jp_{i,j}=p_i$, $\sum_ip_{i,j}=q_j$
 and the min-entropy $H(P)$.

1: { \textbf{set}} $y=(p_1,\cdots,p_{m-1},1,q_1,\cdots,q_{n-1})^T$, Joint-Distr$\leftarrow$ list( ); Joint-Distr-entropy$\leftarrow$ c( ).

2: { \textbf{for}} $k=1,2,\cdots,\binom{mn}{m+n-1}$, { \textbf{do}}

3: $S_k=\{i_1,i_2,\cdots,i_{m+n-1}\}$

4: \textbf{for}~$i=1,\cdots,m$~and~$j=1,\cdots,n$, \textbf{set} $p_{i,j}\leftarrow0$, \textbf{end} .

5: \textbf{for}~$i=1,\cdots,m+n-1$~and~$j=1,\cdots,m+n-1$, \textbf{set} $a_{i,j}\leftarrow0$, \textbf{end}.

6: \textbf{for}~$j=1,\cdots,m+n-1$, \textbf{do}

7: $i_j=(t-1)n+r(1\leq t<m-1, 1\leq r\leq n)$, $a_{t,j}\leftarrow1$, \textbf{end}.

8: \textbf{set} $a_{m,k}\leftarrow1, 1\leq k\leq n$.

9: \textbf{for}~$j=1,\cdots,m+n-1$, \textbf{do}

10: $i_j=(t-1)n+r(1\leq r\leq n-1)$, $a_{m+r,j}\leftarrow1$, \textbf{end}.

11: \textbf{set} $\fA=(a_{i,j})$.

12: \textbf{if}~${\rm det}(\fA)\neq0$ \textbf{then}

13: Solving equation $\fA x=y$, $x=(x_1,\cdots,x_{m+n-1})^T$.

14: \textbf{for}~$j=1,\cdots,m+n-1$, \textbf{do}

15: $i_j=(t-1)n+r$, $p_{t,r}\leftarrow x_j$ ,  \textbf{end}.

16: $P\in $  Joint-Distr, $H(P)\in$~Joint-Distr-entropy.

17: \textbf{end}.

18: {\textbf{find}} minimum value in Joint-Distr-entropy and its corresponding matrix $P$.

19: {\textbf{return}} $(P, H(P))$.

\hrulefill

As examples, some calculating results obtained by the above algorithm will be given in Section 5.

\vskip3mm

\vskip 5mm
\section{Generalizations}
\renewcommand{\theequation}{4.\arabic{equation}}
\setcounter{equation}{0}

In this section, we will generalize the minimum-entropy coupling problem in two directions. Firstly, we study the optimization problem for Schur-concave function on ${\cal C}({\bf p},{\bf q})$. Secondly, we will generalize the minimum-entropy coupling problem to the multi-marginals cases.

\subsection{The optimization problem for Schur-concave function on ${\cal C}({\bf p},{\bf q})$}

From the proof of Theorem~\ref{t1}, one knows that, to obtain the forest structure of a minimal entropy coupling, the strict concave property of the Shannon entropy $H$ is sufficient. In Section 2, a local optimization method is developed, and then, besides the forest structure, other special features, including essential order-preserving and the local order property as revealed in item 2 of Theorem~\ref{t2}, of the minimal entropy coupling are obtained. In the present subsection, we point out that our local optimization method can be generalized to solve the corresponding optimization problem for Schur-concave function on ${\cal C}({\bf p},{\bf q})$.

To introduce the concept of Schur-concave function, we first introduce the concept of {\it majorization}. Note that the concept of majorization plays a key role in constructing proper bounds for the minimum-entropy coupling problem, see \cite{CGV,L} and the references therein.

Recall that for any $x=(x_1,x_2,\ldots,x_m)\in \R^m$, $\bar x=(\bar x_1,\bar x_2,\ldots,\bar x_m)$ be the permutation of $x$ such that $\bar x_1\geq \bar x_2\geq\ldots\geq\bar x_m$, $F_{\bar x}(i), 1\leq i\leq m$ be the the cumulative distribution function defined in (\ref{1.4}).

\bd\label{d4}{\it For any $x,y\in \R^m$, we say $x$ is majorized by $y$, denote by $ x\preceq y$, if
$$F_{\bar x}(i)\leq F_{\bar y}(i), \ {\rm for\ all}\ 1\leq i< m;\ \ F_{\bar x}(m)=F_{\bar y}(m).$$ We say $x$ is strictly majorized by $y$, denote by $x\prec y$, if for some $1\leq i<m$, $F_{\bar x}(i)<F_{\bar y}(i)$.
}\ed
It was proved by Schur \cite{Sch} in 1923 that, $x\preceq y$ if and only if for some {\it doubly stochastic} matrix $D$, \be\label{4.0'} x=Dy.\ee Note that a nonnegative matrix $D$ is called {\it doubly stochastic}, if each row and each column of $D$ sums to unit.

\bd\label{d5} A symmetric function $\Psi: \R^m\rightarrow \R$ is called Schur-convex, if for any $x,y\in\R^m$ with $x\preceq y$, one has $\Psi(x)\leq\Psi(y)$. $\Psi$ is called strict Schur-convex, if for any $x,y\in\R^m$ with $x\prec y$, one has $\Psi(x)<\Psi(y)$. $\Psi$ is called Schur-concave, if $-\Psi$ is Schur-convex.\ed

The Schur-convex property of function is the generalization of the convex property. In fact, by the Birkhoff Theorem \cite{Bir}, the permutation matrices constitute the extreme points of the set of doubly stochastic matrices, note that a permutation matrix is a special matrix obtained from the identity matrix by rearranging rows or columns. That is, if $D$ is doubly stochastic, then there exists permutation matrices $\Pi_i, 1\leq i\leq s$ and $\lambda_i\in (0,1)$ with $\sum_{i=1}^s \lambda_i=1$, such that $D=\sum_{i=1}^s\lambda_i\Pi_i$. Thus, if $x\preceq y$ and $x=Dy$, then for any symmetric convex function $\Psi$, one has $$\Psi(x)=\Psi(Dy)=\Psi\(\(\sum_{i=1}^s\lambda_i\Pi_i\)y\)=\Psi\(\sum_{i=1}^s\lambda_i \(\Pi_iy\)\)\leq \sum_{i=1}^s\lambda_i\Psi(\Pi_iy)=\Psi(y),$$ i.e. $\Psi$ is Schur-convex.

\bt\label{t4} For any ${\bf p}\in{\cal P}^+_m$, ${\bf q}\in {\cal P}_n^+$, suppose $\Psi$ is a strict Schur-concave function on ${\cal C}({\bf p,q})$, then all its minimal value points lie in $\hC({\bf p},{\bf q})$ and
$$\inf_{P\in{\cal C}({\bf p},{\bf q})}\Psi(P)=\min_{P\in \hC({\bf p, q})}H(P).$$ Where $\hC({\bf p,q})={\cal C}_e({\bf p},{\bf q})$ is the extreme point set of ${\cal C}_e({\bf p},{\bf q})$. \et

%


Before giving a proof to Theorem~\ref{t4}, we first give out the following simple version of the local optimization theorem.
\bl\label{l6} For any ${\bf p}\in{\cal P}^+_m$, ${\bf q}\in {\cal P}^+_n$, suppose $\Psi$ is strict Schur-concave on ${\cal C}({\bf p},{\bf q})$. For any $P\in {\cal C}({\bf p},{\bf q})$, let $A=(a_{i,j})_{2\times 2}$ is a 2-nd order submatrix of $P$ satisfying the conditions of Lemma 1, and $A'$ is the 2-nd order matrix obtained from $A$ as in Lemma 1. Let $P'\in {\cal C}({\bf p},{\bf q})$ be the coupling obtained from $P$ by $A'$ taking the place of $A$. Then, as vectors in $\R^{mn}$, one has $P\preceq P'$ and then $\Psi(P')\leq \Psi(P)$. In particular, if $b:=a_{1,2}\wedge a_{2,1}>0$, then $P\prec P'$, $\Psi(P')<\Psi(P)$.  \el

{\it Proof.} Without loss of generality, assume $a_{1,1}\geq a_{2,1}\geq a_{1,2}$, $a_{1,1}\geq a_{2,2}$. Note that in this case, one has $b=a_{1,2}$, $a'_{1,2}=0$ and
$$A=(a_{i,j})=\(\ba{ll}a_{1,1} &a_{1,2}\\[5mm]
a_{2,1} & a_{2,2}\ea\),\ A'=(a'_{i,j})=\(\ba{ll}a_{1,1}+b &0\\[5mm]
a_{2,1}-b & a_{2,2}+b\ea\).$$
Denote by $x=(a_{1,1},a_{2,1},a_{1,2},a_{2,2})^T$, $y=(a_{1,1}+b, a_{2,1}-b,a_{2,2}+b,0)^T$. We claim that $x\preceq y$. Actually, it always holds that $F_{\bar x}(1)\leq F_{\bar y}(1)$, $F_{\bar x}(3)\leq F_{\bar y}(3)$ and $F_{\bar x}(4)=F_{\bar y}(4)$, it only remains to prove $F_{\bar x}(2)\leq F_{\bar y}(2)$. In the case of $a_{2,2}\leq a_{2,1}$, one has $F_{\bar x}(2)=a_{1,1}+a_{2,1}\leq (a_{1,1}+b+a_{2,1}-b)\vee (a_{1,1}+b+a_{2,2}+b)=F_{\bar y}(2)$; in the case of $a_{1,1}\geq a_{2,2}\geq {a_{2,1}}$, one has $F_{\bar x}(2)=a_{1,1}+a_{2,2}\leq a_{1,1}+b+a_{2,2}+b=F_{\bar y}(2)$. If $b>0$, then $F_{\bar x}(1)=a_{1,1}< F_{\bar y}(1)=a_{1,1}+b$, and $x\prec y$.

Now, let $D$ be the doubly stochastic matrix such that $x=Dy$, let $\Pi$ be the $mn$-order permutation matrix such that, as vectors in $\R^{mn}$, $\Pi P=(x^T,z)^T$ and $\Pi P'=(y^T,z)^T$, where $z=(z_1,z_2,\ldots,z_{mn-4})\in \R^{mn-4}$. Let
$$\bar D=\(\ba{ll}D & 0\\[2mm]
0 & I\ea\),$$ where $I$ is the $(mn-4)$-order identity matrix. Then $\bar D$ is doubly stochastic and $$\Pi P=(x^T,z)^T=\bar D(y^T,z)^T=\bar D\Pi P',$$ thus $P=\Pi^{-1}\bar D \Pi P'$. Since $\Pi^{-1}\bar D\Pi$ is doubly stochastic, it follows from (\ref{4.0'}) that $P\preceq P'$. Clearly, if $b>0$, then $P\prec P'$.\QED

%
%


{\it Proof of Theorem~\ref{t4}:} Suppose $\Psi$ is a strict Schur-concave function on ${\cal C}({\bf p},{\bf q})$ and $\tilde P\in {\cal C}({\bf p},{\bf q})$ is a minimal value point. Then by Lemma~\ref{l6} and the same argument in the proof of Theorem~\ref{t2}, $V(\tilde P)$ is a complete forest. By Proposition~\ref{p0}, i), there exists some $T\in{\cal T}$ such that $V(\tilde P)\subset T$. Namely, $\tilde P$ is consistent to $T$ and $\tilde P\in {\hC}({\bf p}, {\bf q})$, thus we finish the proof.\QED

At the end of this subsection, we introduce the concept of $(\Phi,\hbar)$-entropy, which consists a large class of strict Schur-concave functions including the Shannon entropy, the R\'enyi entropy and the Tsallis entropy.

Let $\Phi: \R \longrightarrow\R$, $\hbar:[0,1]\longrightarrow \R$ are two functions. For any $0\leq c\leq 1$, define
\be\label{4.1}\hbar_c(x):=\hbar(x)+\hbar(c-x),\ x\in[0,c].\ee
In this subsection, we will consider the function pairs $(\Phi,\hbar)$ satisfying the following monotonicity conditions:
\be\label{4.0}
\ba{ll}
\hskip -5.5cm \bullet & \Phi {\rm \ is\ strictly\ monotone\ }; \\[2mm]
\hskip -5.5cm \bullet & {\rm for\ any\ } 0\leq c\leq 1, \hbar_c{\rm \ is\ strictly\ monotone\ in\ } [0,c/2];\\[2mm]
\hskip -5.5cm \bullet & {\rm for\ any\ } 0\leq c\leq 1, \Phi(\hbar_c){\rm \ is\ strictly\ increasing\ in\ } [0,c/2];\\[2mm]
\ea
\ee
\bd\label{d3} Suppose $(\Phi,\hbar)$ is a function pair satisfying (\ref{4.0}). For any ${\bf p}\in {\cal P}_m$, one kind of entropy of $\bf p$, denote by $H(\bf p)$, is called a $(\Phi,\hbar)$-entropy, if $H(\bf p)$ can be written as
\be\label{4.2}H({\bf p})=\Phi\(\sum_{i=1}^m\hbar(p_i)\).\ee For any ${\bf p}\in {\cal P}_m$, ${\bf q}\in{\cal P}_n$, $P\in {\cal C}({\bf p, q})$, the $(\Phi,\hbar)$-entropy of $P$ is given by
\be H(P)=\Phi\(\sum_{i=1}^m\sum_{j=1}^n\hbar(p_{i,j})\).\ee\ed

Clearly, the Shannon entropy is a $(\Phi,\hbar)$-entropy with $\Phi(x)=x$, $\hbar(x)=-x\log x$. Furthermore, for $\a\geq 0,\ \a\not=1$, the R\'enyi entropy defined in (\ref{4.5})
is the $(\Phi,\hbar)$-entropy with $\Phi(x)={\log x}/{(1-\a)}$, $\hbar(x)=x^\a$; the Tsallis entropy defined in (\ref{4.6})
is the $(\Phi,\hbar)$-entropy with $\Phi(x)=x/{(1-\a)}$, $\hbar(x)=x^\a-x$.

\bp {\it Suppose $H$ is a $(\Phi,\hbar)$-entropy with differentiable function pair $(\Phi,\hbar)$, then for any ${\bf p}\in {\cal P}^+_m$, ${\bf q}\in{\cal P}^+_n$, $H$ is strict Schur-concave on ${\cal C}({\bf p},{\bf q})$.}\ep

{\it Proof.} For a symmetric differentiable function $\Psi:\R^m\rightarrow \R$, $\Psi$ is Schur-concave, if and only if the following Schur-Ostrowski condition \cite{Os} holds:
$$(x_i-x_j)\(\frac{\partial \Psi}{\partial x_i}-\frac{\partial \Psi}{\partial x_j}\)\leq 0,\ {\rm for\ any}\ 1\leq i\not=j\leq m.$$ Now we have
$H(x)=\Phi\(\sum_{i=1}^m\hbar(x_i)\)$ and then, by the monotonicity condition (\ref{4.0}),
$$\ba{ll}(x_i-x_i)\Dp\(\frac{\partial H}{\partial x_i}-\frac{\partial H}{\partial x_j}\)&=\Dp(x_i-x_j)\Phi'\(\sum_{i=1}^m\hbar(x_i)\)\(\hbar'(x_i)-\hbar'(x_j)\)\\[4mm]
&=\Dp(x_i-x_j)\Phi'\(\sum_{i=1}^m\hbar(x_i)\)\hbar_{x_i+x_j}'(x_i)\leq 0,\ea$$ for all $1\leq i\not=j\leq m$ and $x=(x_1,x_2,\ldots,x_m)\in \R^m_+$. The above inequality holds strictly if $x_i\not=x_j$. Thus we finish the proof.\QED

We finish the subsection by giving an example to show that a $(\Phi,\hbar)$-entropy $H$ can be not concave. To this end, let $\hbar':[0,1]\rightarrow (0,\infty)$ be the continuously differentiable function such that
\be\label{4.e1}\hbar'(x)\left\{\ba{ll}=1-x,&{\rm if}\ 0\leq x\leq 5/8;\\[2mm]
<1/2,&{\rm if}\ 5/8\leq x\leq 7/8;\\[2mm]
=2x-13/8,&{\rm if}\ 7/8\leq x\leq 1.\ea \right.\ee Define $\hbar(x):=\int_{0}^x\hbar'(y)dy,\ x\in[0,1]$, and let $\Phi(x)=x, \ x\in \R$.

For any $c\in [0,1]$, let $\hbar_c(x)=\hbar(x)+\hbar(c-x)$, $x\in[0,c]$. Then, for any $x\in[0,c/2]$, by the definition of $\hbar'$, one has $\hbar'_c(x)=\hbar'(x)-\hbar'(c-x)>0$. Thus $(\Phi,\hbar)$ is a function pair satisfying the monotonicity condition (\ref{4.0}), and the $(\Phi,\hbar)$-entropy $H$ is well defined by Definition~\ref{d3}.

Let's consider the $(\Phi,\hbar)$-entropy $H$ on ${\cal P}_2$. For any ${\bf p}=(p_1,p_2)\in {\cal P}_2$, without loss of generality, suppose that $p_1\leq p_2$. According to Definition~\ref{d3},
\be\label{4.e2} H({\bf p})=\Phi\(\sum_{i=1}^2\hbar(p_i)\)=\hbar_1(p_1).\ee

However, for any $x\in[0,1/2]$, by (\ref{4.e1}), one has
$$\hbar_1''(x)=\hbar''(x)+\hbar''(1-x)=\left\{\ba{ll}-1+2=1,&{\rm if}\ 0\leq x\leq 1/8;\\[2mm]
-1+-1=-2,&{\rm if}\ 3/8\leq x\leq 1/2.\ea\right.$$
Thus, the $(\Phi,\hbar)$-entropy $H$ defined in (\ref{4.e2}) is not a concave function on ${\cal P}_2$.

Finally, for $m,n\geq 2$, ${\bf p}\in {\cal P}^+_m, {\bf q}\in {\cal P}^+_n$. If $\max\{p_1,\ldots,p_m,q_1,\ldots,q_n\}\leq 5/8$, then by (\ref{4.e1}), the $(\Phi,\hbar)$-entropy $H$ defined by
$$H(P)=\Phi\(\sum_{i=1}^m\sum_{j=1}^n\hbar(p_{i,j})\)=\sum_{i=1}^m\sum_{j=1}^n\hbar(p_{i,j})$$ is strict concave on ${\cal C}({\bf p},{\bf q})$.


\subsection{The minimum-entropy coupling problem for multi-marginal cases}

The minimum-entropy coupling problem (\ref{1}) has been naturally generalized to the following multi-marginal case by mathematicians.

For any integer $d\geq 2$, for any integers $m_1,m_2,\ldots,m_d\geq 2$, and for any probability distributions ${\bf p}^1\in{\cal P}^+_{m_1}$, ${\bf p}^2\in{\cal P}^+_{m_2}$, $\ldots$ , ${\bf p}^d\in{\cal P}^+_{m_d}$, write $\mathbb S=\{{\bf p}^i:1\leq i\leq d\}$ and denote by ${\cal C}(\mathbb S)$ the collection of couplings of $\{{\bf p}^i:1\leq i\leq d\}$, denote by ${\cal C}_e(\mathbb S)$ the set of extreme points of ${\cal C}(\mathbb S)$. For any $P=(p_{l_1,l_2,\ldots, l_d})_{m_1\times m_2\times\ldots\times m_d}\in {\cal C}(\mathbb S)$, let $H(P)$ be the $(\Phi,\hbar)$-entropy of $P$ given in Definition~\ref{d3}, i.e.
$$H(P)=\Phi\(\sum_{l_1=1}^{m_1}\sum_{l_2=1}^{m_2}\cdots\sum_{l_d=1}^{m_d}\hbar(p_{l_1,l_2,\ldots, l_d})\).$$
Then the minimum-entropy coupling problem for marginals ${\bf p}^1,{\bf p}^2,\ldots,{\bf p^d}$ is the following optimization problem:
\be\label{4.7} \tilde P:H(\tilde P)=\inf_{P\in{\cal C}(\mathbb S)}H(P).\ee
Here we declare that the solving procedure for the above optimization problem (\ref{4.7}) is completely similar to that of problem (\ref{1}). In the rest of this subsection, we only state the results and omit the detailed proofs.

Let $G_{m_1,\ldots,m_d}$ be the graph with vertex set $V_{m_1,\ldots,m_d}=[m_1]\times\ldots\times[m_d]\subset \Z^d$, the $d$-dimensional integer lattice,  and edge set $E_{m_1,\ldots,m_d}$, a collection of edge $e=\langle u,v\rangle$ such that $u$ only differs from $v$ at one coordinate. Completely similar to Definition~\ref{d1}, we define \textbf{continuous}, \textbf{directed} path and \textbf{circuit} in $G_{m_1,\ldots,m_d}$, and for subgraph of $G_{m_1,\ldots,m_d}$, we introduce the concepts of \textbf{forest}, \textbf{tree} and \textbf{completeness}. Then, similar to Proposition~\ref{p0}, we have

\bp\label{p3}Suppose $V\subset V_{m_1,\ldots,m_d}$, then
\begin{description}
  \item[\ \ \ \ \it i)] if $V$ is a forest with $k$ connected components, then $V$ is complete if and only if $|V|=\Dp\sum_{i=1}^d m_i-d-k+2$;
  \item[\ \ \ \ \it ii)] if $|V|=\Dp\sum_{i=1}^d m_i-(d-1)$, then $V$ is a tree if and only if $V$ is complete and connected;
  \item[\ \ \ \ \it iii)] if $V$ is complete and connected, then $|V|\geq \Dp\sum_{i=1}^d m_i-(d-1)$ and $|V|=\Dp\sum_{i=1}^d m_i-(d-1)$ if and only if $V$ is a tree.
\end{description}
\ep

Denote ${\cal T}_{m_1,\ldots,m_d}=\{T\subset V_{m_1,\ldots,m_d}:$ $T$ is a tree and $|T|=\sum_{i=1}^d m_i-(d-1)\}$. For any $T\in{\cal T}_{m_1,\ldots,m_d}$ and $P\in{\cal C}(\mathbb S)$, we call $P$ is \textbf{consistent} with $T$, if $V(P)\subset T$, where $V(P)=\{(i_1,i_2,\ldots,i_d)\in V_{m_1,\ldots,m_d}:p_{i_1,i_2,\ldots,i_d}>0\}$ is the support of $P$.

\bp\label{p4} For any $T\in{\cal T}_{m_1,\ldots,m_d}$, there exists at most one $P\in{\cal C}(\mathbb S)$, such that $P$ is consistent with $T$.\ep
Now, let ${\mathscr C}(\mathbb S)=\{P\in{\cal C}(\mathbb S):$ for some $T\in{\cal T}_{m_1,\ldots,m_d}$, $P$ is consistent with $T$$\}$. Clearly $$|{\mathscr C}(\mathbb S)|\leq\binom{\prod_{i=1}^dm_i}{\sum_{i=1}^dm_i-(d-1)}.$$

\bt\label{t8}For any $d\geq 2$, for any $m_1,m_2,\ldots,m_d\geq 2$, and for any probability distributions ${\bf p}^1\in{\cal P}^+_{m_1}$, ${\bf p}^2\in{\cal P}^+_{m_2}$, $\ldots$ , ${\bf p}^d\in{\cal P}^+_{m_d}$. Then
\be {\cal C}_e(\mathbb S)={\hC}(\mathbb S).\ee
If $\tilde P$ solves the optimization problem (\ref{4.7}), then $\tilde P\in {\mathscr C}(\mathbb S)$ and
\be\label{4.8}H(\tilde P)=\min_{P\in{\mathscr C}(\mathbb S)}H(P).\ee
\et

Theorem~\ref{t8} can be proved in two steps. Step 1, by updating the local optimization theorem (Theorem~3 and Lemma~\ref{l6}) to a general version, we prove that, for any minimum-entropy coupling $P\in{\cal C}(\mathbb S)$, $V(P)$ is a forest, then for some tree $T\in{\cal T}_{m_1,\ldots,m_d}$, $P$ is consistent with $T$ and hence $P\in{\mathscr C}(\mathbb S)$. Let $e_i=(0,\ldots,0,1,0,\ldots,0)$, $i=1,2,\ldots,d$, be the $i$-th coordinate unit vector in $\R^d$, for any $2\leq k\leq d-1$, for $1\leq i_1<i_2<\dots<i_k\leq d$, denote by $Hyp(i_1,i_2,\ldots,i_k)$ the $k$-dimensional coordinate hyperplane spanned by vector family $\{e_{i_j}:j=1,2,\ldots,k\}$. Now, for any $P\in {\cal C}(\mathbb S)$, suppose $A$ is a $2\times 2$ submatrix of $P$, which lies in a $2$-dimensional hyperplane parallel to some coordinate hyperplane $Hyp(i,j)$, $1\leq i<j\leq d$. Let $P'$ be the matrix obtained from $P$ by $A'$ taking the place of $A$, where $A'$ is obtained from $A$ as in Lemmas 1, 2 and~\ref{l1}. To update the local optimization theorem to the general case, it suffices to show the fact that $P'\in{\cal C}(\mathbb S)$. To this end, let's consider the $(d-1)$-dimensional hyperplane in $V_{m_1,\ldots,m_d}$:
\be\label{4.9} Hyp(t,z):=\{(l_1,l_2,\ldots,l_d)\in V_{m_1,\ldots,m_d}:l_t=z\},\ 1\leq t\leq d,\ 1\leq z\leq m_t.\ee Clearly $Hyp(t,z)$ is parallel to the $(d-1)$-dimensional coordinate hyperplane $Hyp(1,\ldots,t-1,t+1,\ldots,d)$, and
$$\sum_{(l_1,\ldots,\l_d)\in H(t,z)}p_{l_1,\ldots,l_d}=p^t_z,$$ the $z$-th component of distribution ${\bf p}^t$. Now, let's consider the possible relative position between submatrix $A$ and the $(d-1)$-dimensional hyperplane $H(t,z)$. In the case of $i,j\not=t$, either all entries of $A$ lie in $Hyp(t,z)$ or no entry of $A$ lies in $Hyp(t,z)$; in the case of $i=t$ (resp. $j=t$), either only  two entries of $A$, which lie in a line parallel to vector $e_j$ (resp. $e_i$), lie in $Hyp(t,z)$ or no entry of $A$ lies in $Hyp(t,z)$. Then by the definition of $P'$, one always has
$$\sum_{(l_1,\ldots,\l_d)\in Hyp(t,z)}p'_{l_1,\ldots,l_d}=\sum_{(l_1,\ldots,\l_d)\in Hyp(t,z)}p_{l_1,\ldots,l_d}=p^t_z,$$ where $p'_{l_1,\ldots,l_d}$ is the entry in $P'$. Thus we obtain $P'\in{\cal C}(\mathbb S)$.

Step 2, by proving Proposition~\ref{p4}, we obtain $|{\mathscr C}(\mathbb S)|<\infty$. The proof of Proposition~\ref{p4} is similar to that of Proposition~\ref{p-}, but is more complicated. For any $P\in {\cal C}(\mathbb S)$ such that for some $T\in{\cal T}_{m_1,\ldots,m_d}$, $V(P)\subset T$, suppose vertex $(l_1,l_2,\ldots,l_d)$ is a leaf of $T$ ($T$ has at least two leaves unless $|T|=1$). The key fact for a proof to Proposition~\ref{p4} is that, $p_{l_1,l_2,\ldots,l_d}$, the entry of $P$, is completely determined by $T$ and the set of marginals $\mathbb S$. In fact, since $(l_1,l_2,\ldots,l_d)$ is a leaf of $T$, then $p_{l_1,\ldots,l_d}$ is the unique nonnegative element in some $(d-1)$-dimensional hyperplane of $V_{m_1,\ldots,m_d}$ as given in (\ref{4.9}). Without loss of generality, suppose this hyperplane is parallel to $Hyp(1,2,\ldots,t-1,t+1,\ldots,d)$, then $p_{l_1,l_2,\ldots,l_d}=p^t_{l_t}$. Let $G'=G_{m_1,\ldots,m_{t-1},m_t-1,m_{t+1},\ldots,m_d}$ be the graph obtained from $G_{m_1,\ldots,m_d}$ by deleting all vertices in this $(d-1)$-dimensional hyperplane and all relevant edges, let $T'=T\setminus\{(l_1,l_2,\ldots,l_d)\}$, then $T'$ is a tree in $G'$. Repeat the above procedure $|T|$ times, all entries of $P$ are determined.

Finally, for any $V\subset V_{m_1,\ldots,m_d}$ with $|V|=\sum_{i=1}^dm_i-(d-1)$, similar to Definition~\ref{d10}, one can define the structure matrix ${\cal A}(V)$ such that $V$ is a tree if and only if $\det({\cal A}(V))\not=0$. Then, a similar but more complicated algorithm follows. In the next section, we will give some calculating results for $d=3$ and small $m_1,m_2,m_3$.

\vskip 5mm
\section{Examples}
\renewcommand{\theequation}{5.\arabic{equation}}
\setcounter{equation}{0}

In this section, as examples, by using Theorem~\ref{t4'} and the algorithm given in Section 3, we first give out
some calculating results for the classical minimum-entropy coupling problem (\ref{1}) for $m,n\leq 5$. For the problem is essentially NP-hard, unfortunately, we can not obtain a result for $m,n\geq 6$ by using a personal computer. Note that in all these examples, we choose $2$ as the base of the $\log$-function.
\vskip 2mm
The following Examples \ref{e5.1}-\ref{e5.5} are calculating results for Shannon entropy.

\bex\label{e5.1} {Case m=n=3}:
\vskip 3mm
1, if ${\bf p}=(0.50, 0.40,0.10),\ {\bf q}=(0.60, 0.20, 0.20)$, then
$$
\tilde P=\begin{pmatrix}
0.50&0&0\\
0&0.20&0.20\\
0.10&0&0
\end{pmatrix},
\ \  H(\tilde P)=1.760964.$$

2, if ${\bf p}=(0.40, 0.35,0.25),\ {\bf q}=(0.38, 0.34, 0.28)$, then
$$
\tilde P=\begin{pmatrix}
0.38&0&0.02\\
0&0.34&0.01\\
0&0&0.25
\end{pmatrix},
\ \  H(\tilde P)=1.738942.$$\eex


\bex\label{e5.2}{Case m=n=4}:
\vskip 3mm
1, if ${\bf p}=(0.40, 0.30,0.20,0.10),\ {\bf q}=(0.38, 0.27, 0.20, 0.15)$, then
$$
\tilde P=\begin{pmatrix}
0.38&0&0&0.02\\
0&0.27&0&0.03\\
0&0&0.20&0\\
0&0&0&0.10
\end{pmatrix},
\ \  H(\tilde P)=2.101697.$$

2, if ${\bf p}=(0.50, 0.20,0.18,0.12),\ {\bf q}=(0.45, 0.25, 0.16, 0.14)$, then
$$
\tilde P=\begin{pmatrix}
0.45&0.05&0&0\\
0&0.20&0&0\\
0&0&0.16&0.02\\
0&0&0&0.12
\end{pmatrix},
\ \  H(\tilde P)=2.101845.$$\eex



\bex\label{e5.3}{Case m=5, n=4}:
\vskip 3mm
1, if ${\bf p}=(0.43, 0.30,0.15,0.10,0.02),\ {\bf q}=(0.40, 0.30, 0.18, 0.12)$, then
$$
\tilde P=\begin{pmatrix}
0.40&0&0.03&0\\
0&0.30&0&0\\
0&0&0.15&0\\
0&0&0&0.10\\
0&0&0&0.02
\end{pmatrix},
\ \ H(\tilde P)=2.057242.$$

2, if ${\bf p}=(0.70, 0.15,0.10,0.03,0.02),\ {\bf q}=(0.50, 0.20, 0.17, 0.13)$, then
$$
\tilde P=\begin{pmatrix}
0.50&0.20&0&0\\
0&0&0.15&0\\
0&0&0&0.10\\
0&0&0&0.03\\
0&0&0.02&0
\end{pmatrix},
\ \ H(\tilde P)=1.971767.$$\eex

\bex\label{e5.4}{Case m=n=5}:
\vskip 3mm
1, if ${\bf p}=(0.33, 0.22,0.17,0.16,0.12),\ {\bf q}=(0.30, 0.25, 0.20, 0.15,0.10)$, then
$$
\tilde P=\begin{pmatrix}
0.30&0.03&0&0&0\\
0&0.22&0&0&0\\
0&0&0.17&0&0\\
0&0&0.01&0.15&0\\
0&0&0.02&0&0.10
\end{pmatrix},
\ \  H(\tilde P)=2.51007.$$

2, if ${\bf p}=(0.40, 0.30,0.15,0.10,0.05),\ {\bf q}=(0.28, 0.27, 0.21, 0.16,0.08)$, then
\be\label{5.1}
\tilde P=\begin{pmatrix}
0.28&0&0.12&0&0\\
0&0.27&0&0&0.03\\
0&0&0&0.15&0\\
0&0&0.09&0.01&0\\
0&0&0&0&0.05
\end{pmatrix},
\ \  H(\tilde P)=2.54881.\ee
\eex

%

\bex\label{e5.5} In case $m=n=5$, we give two examples to reveal the non-uniqueness of the minimal entropy coupling:
\vskip 3mm
1, if ${\bf p}=(0.50, 0.30,0.08,0.07,0.05), {\bf q}=(0.35, 0.25, 0.20, 0.14,0.06)$, then
\be\label{5.2}
\tilde P=\begin{pmatrix}
0.35&0&0.15&0&0\\
0&0.25&0.05&0&0\\
0&0&0&0.08&0\\
0&0&0&0.06&0.01\\
0&0&0&0&0.05
\end{pmatrix}\ {\rm or\ }\begin{pmatrix}
0.35&0&0.15&0&0\\
0&0.25&0&0&0.05\\
0&0&0&0.08&0\\
0&0&0&0.06&0.01\\
0&0&0.05&0&0
\end{pmatrix},\ee
 $H(\tilde P)=2.474319.$

2, if ${\bf p}=(0.55, 0.35,0.05,0.03,0.02), {\bf q}=(0.40, 0.30, 0.20, 0.06,0.04)$, then
\be\label{5.3}
\tilde P=\begin{pmatrix}
0.40&0&0.15&0&0\\
0&0.30&0.05&0&0\\
0&0&0&0.05&0\\
0&0&0&0&0.03\\
0&0&0&0.01&0.01
\end{pmatrix}\ {\rm or\ }\begin{pmatrix}
0.40&0&0.15&0&0\\
0&0.30&0&0.05&0\\
0&0&0.05&0&0\\
0&0&0&0&0.03\\
0&0&0&0.01&0.01
\end{pmatrix},\ee
$H(\tilde P)=2.177242.$\eex

\br It seems that, in most cases, for the minimum-entropy coupling $\tilde P$, $\tilde P(\sigma,\pi)$ may be $\kappa({\bf p,q})$-blocked as in (\ref{2.6}) for some $(\sigma,\pi)\in\Sigma_m\times\Sigma_n$, see the above calculating results obtained in Example~\ref{e5.1}, 1, Example \ref{e5.2}, Example \ref{e5.3}, 2, Example \ref{e5.4}, 1 and Example~\ref{e5.5}. Of course, this is not always true, Example~\ref{e5.4}, 2 is a counter example.\er

The following Examples \ref{e5.7} and \ref{e5.8} are calculating results for R\'enyi entropy and Tsallis entropy. Here we denote $H^R_\a (P)$, $H^T_\a (P)$ the R\'enyi entropy, the Tsallis entropy (with parameter $\a$) of $P$ respectively.

\bex\label{e5.7} In the case $m=n=5$, for $\a=0.1,\, 0.5,\, 0.9,\, 1.1,\, 1.5$ and $2.0$, we calculate the corresponding minimal joint entropies respectively.

1, if ${\bf p}=(0.50, 0.30,0.08,0.07,0.05)$, ${\bf q}=(0.35, 0.25, 0.20,$$ 0.14,0.06)$, i.e. the same ${\bf p}, {\bf q}$ as in Example~\ref{e5.5}, 1, then the minimum-entropy couplings $\tilde P$'s are the same as given in (\ref{5.2}) and the corresponding entropy values are given in the following table.
\vskip 2mm
\begin{center}
\begin{tabular}{|c|c|c|c|c|c|c|}
  \hline
  $\a$& 0.1 & 0.5 & 0.9 & 1.1 & 1.5 & 2.0\\
  \hline
  $H^R_\a (\tilde P)$& 2.935792 & 2.705417 & 2.515795 & 2.435067 & 2.298609 & 2.167475 \\
  $H^T_\a (\tilde P)$& 5.796255 & 3.107823 & 1.905098 & 1.553103 & 1.098315 & 0.7774 \\
  \hline
\end{tabular}
\end{center}

2, if ${\bf p}=(0.55, 0.35,0.05,0.03,0.02)$, ${\bf q}=(0.40, 0.30, 0.20, 0.06,0.04)$, i.e. the same ${\bf p}, {\bf q}$ as in Example~\ref{e5.5}, 2, then the minimum-entropy couplings $\tilde P$'s are the same as given in (\ref{5.3}) and the corresponding entropy values are given in the following table.
\vskip 2mm
\begin{center}
\begin{tabular}{|c|c|c|c|c|c|c|}
  \hline
  $\a$& 0.1 & 0.5 & 0.9 & 1.1 & 1.5 & 2.0\\
  \hline
  $H^R_\a (\tilde P)$& 2.891993 & 2.511479 & 2.232101 & 2.127567 & 1.971572 & 1.843733\\
  $H^T_\a (\tilde P)$& 5.638465 & 2.77579 & 1.673281 & 1.371132 & 0.9900989 & 0.7214 \\
  \hline
\end{tabular}
\end{center}
\eex

\bex\label{e5.8} For ${\bf p}=(0.40, 0.30,0.15,0.10,0.05)$, ${\bf q}=(0.28, 0.27, 0.21, 0.16,0.08)$, i.e. the same ${\bf p},{\bf q}$ as in Example~\ref{e5.4}, 2,
for $\a=0.1,\,0.5,\,0.9,\,1.1$ and $1.5$, the minimum-entropy coupling $\tilde P$ is the same as given in (\ref{5.1}) and the corresponding entropy values are given in the following table.
\vskip 2mm
\begin{center}
\begin{tabular}{|c|c|c|c|c|c|}
  \hline
  $\a$& 0.1 & 0.5 & 0.9 & 1.1 & 1.5\\
  \hline
  $H^R_\a (\tilde P)$& 2.93921 & 2.733940 & 2.580667 & 2.519114 & 2.418465 \\
  $H^T_\a (\tilde P)$& 5.840234 & 3.158572 & 1.958751 & 1.602169 & 1.135003 \\
  \hline
\end{tabular}
\end{center}
But for $\a=2.0$, the minimum-entropy coupling is
\be\label{5.4}
\tilde P=\begin{pmatrix}
0&0.27&0.13&0&0\\
0.28&0&0.01&0.01&0\\
0&0&0&0.15&0\\
0&0&0.02&0&0.08\\
0&0&0.05&0&0
\end{pmatrix}
\ee with $H^R_{2.0}(\tilde P)=2.320486$, $H^T_{2.0}(\tilde P)=0.7998$.
\eex

At the end of this section, we give out some calculating results for the minimum-(Shannon) entropy coupling in multi-marginal cases. In the following Examples~\ref{e5.9} and \ref{e5.10}, we choose $d=3$, $m_1=m_2=m_3=3$.

\bex\label{e5.9} For ${\bf p}=(0.50,0.40,0.10)$, ${\bf q}=(0.60,0.20,0.20)$ and ${\bf r}=(0.40,0.30,0.30)$, the minimum-entropy coupling is $\tilde P=(p_{i,j,r})_{3\times3\times 3}$ with
$$(p_{i,j,1})=\begin{pmatrix}
0&0&0\\
0.40&0&0\\
0&0&0
\end{pmatrix};\,(p_{i,j,2})=\begin{pmatrix}
0.10&0.20&0\\
0&0&0\\
0&0&0
\end{pmatrix};\,(p_{i,j,3})=\begin{pmatrix}
0&0&0.20\\
0&0&0\\
0.10&0&0
\end{pmatrix}
$$
and $H(\tilde P)=2.121928$. Note that ${\bf p},{\bf q}$ are the same as in Example \ref{e5.1}, 1, the marginal coupling of ${\bf p}$ and ${\bf q}$ in $\tilde P$ is
$$\begin{pmatrix}
0.10&0.20&0.20\\
0.40&0&0\\
0.10&0&0
\end{pmatrix},$$which \textbf{differs} from the optimal coupling given in Example \ref{e5.1}, 1.
\eex

\bex\label{e5.10} For ${\bf p}=(0.40,0.35,0.25)$, ${\bf q}=(0.38,0.34,0.28)$ and ${\bf r}=(0.45,0.35,0.20)$, the minimum-entropy coupling is $\tilde P=(p_{i,j,r})_{3\times3\times 3}$ with
$$(p_{i,j,1})=\begin{pmatrix}
0.38&0&0.02\\
0&0&0\\
0&0&0.05
\end{pmatrix};\,(p_{i,j,2})=\begin{pmatrix}
0&0&0\\
0&0.34&0.01\\
0&0&0
\end{pmatrix};\,(p_{i,j,3})=\begin{pmatrix}
0&0&0\\
0&0&0\\
0&0&0.20
\end{pmatrix}
$$
and $H(\tilde P)=1.919424$. Here ${\bf p},{\bf q}$ are the same as in Example \ref{e5.1}, 2, the marginal coupling of ${\bf p}$ and ${\bf q}$ in $\tilde P$ is
$$\begin{pmatrix}
0.38&0&0.02\\
0&0.34&0.01\\
0&0&0.25
\end{pmatrix},$$which \textbf{coincides} with the optimal coupling given in Example \ref{e5.1}, 2.
\eex

The following Example~\ref{e5.11} is a calculating result for $d=3$, $m_1=2$, $m_2=3$ and $m_3=4$.
\bex\label{e5.11} For ${\bf p}=(0.30,0.70)$, ${\bf q}=(0.10,0.40,0.50)$ and ${\bf r}=(0.15,0.20,0.25,0.40)$, the minimum-entropy coupling is $\tilde P=(p_{i,j,r})_{2\times3\times 4}$ with
$$(p_{i,j,1})=\begin{pmatrix}
0&0&0.05\\
0.10&0&0
\end{pmatrix};\,(p_{i,j,2})=\begin{pmatrix}
0&0&0\\
0&0&0.20
\end{pmatrix};\,(p_{i,j,3})=\begin{pmatrix}
0&0&0.25\\
0&0&0
\end{pmatrix};\,(p_{i,j,4})=\begin{pmatrix}
0&0&0\\
0&0.40&0
\end{pmatrix}
$$
and $H(\tilde P)=2.041446$.\eex


\end{document}